\RequirePackage{fix-cm}
\documentclass[smallextended]{svjour3}       
\smartqed  
\usepackage{graphicx}
\usepackage{amssymb}
\usepackage{graphicx,fancyhdr}
\usepackage{multirow,amsmath}
\usepackage{dsfont}
\usepackage{subfigure}
\usepackage{bbm}
\usepackage{booktabs}
\usepackage{bbm}
\usepackage{amsfonts}
\usepackage{mathrsfs,float}
\usepackage{color,bm}

\begin{document}

\title{Fourier Extension  Based on Weighted Generalized Inverse
}


\author{Zhenyu Zhao  \and Yanfei Wang \and Anatoly G. Yagola\and Xusheng Li
}


\institute{Zhenyu Zhao \at
              School of Mathematics and Statistics, Shandong University of Technology, Zibo, 255049, China \\
                          \email{wozitianshanglai@163.com}           
           \and
          Yanfei Wang \at
              Key Laboratory of Deep Petroleum Intelligent Exploration and Development,
Institute of Geology and Geophysics, Chinese Academy of Sciences, Beijing, 100029,China\\
              \email{yfwang@mail.iggcas.ac.cn}
              \and
              Anatoly G. Yagola \at
              Department of Mathematics, Faculty of Physics, Lomonosov Moscow State
University, Vorobyevy Gory, 119991 Moscow, Russia\\
  \email{yagola@physics.msu.ru}
  \and
  Xusheng Li \at School of Mathematics and Statistics, Shandong University of Technology, Zibo, 255049, China\\
    \email{lixusheng1997@163.com}
}

\date{Received: date / Accepted: date}

\maketitle

\begin{abstract}
This paper introduces a weighted generalized inverse framework for Fourier extensions, designed to suppress spurious oscillations in the extended region while maintaining high approximation accuracy on the original interval. By formulating the Fourier extension problem as a compact operator equation, we propose a weighted best-approximation solution that incorporates a priori smoothness information through suitable weight operators on the Fourier coefficients. This leads to a regularization scheme based on the generalized truncated singular value decomposition (GTSVD). Under algebraic and exponential smoothness assumptions, convergence analysis demonstrates optimal $L^2$ accuracy and improved stability for derivatives. Compared with classical Fourier extension using standard TSVD, the proposed method effectively controls high-frequency components and yields smoother extensions. A practical discretization using uniform sampling is developed, along with an adaptive design of weight functions. Numerical experiments confirm that the method significantly improves derivative approximations and reduces oscillations in the extended domain without compromising accuracy on the original interval.
\end{abstract}

\section{Introduction}
Due to its flexibility, frame approximation has received increasing attention in recent years. Its related theories and numerical algorithms are being widely studied and have made substantial progress \cite{adcock2019frames,adcock2020frames,adcock2020approximating,adcock2021frame,herremans2024efficient,ebner2023regularization}.
Fourier extension is an important form of frame approximation, which has significantly contributed to the advancement of the field of frame theory. Initially proposed to mitigate the Gibbs phenomenon observed in classical Fourier series, Fourier extension has since witnessed extensive research. In recent years, many scholars have made substantial contributions to developing its mathematical foundations, devising efficient algorithms, reduce artifacts, enhance the resolution of inversion results and exploring its applications across various fields \cite{Boyd2002A,Huybrechs2010,Lyon2011,Matthysen2016FAST,2012Sobolev,Matthysen2018,Bruno2022,Webb2020}.

Similar to general frame approximations, Fourier extension also encounters challenges stemming from its near-redundancy, which typically manifests as severe ill-conditioning in the least-squares system for computing the approximation coefficients. Nevertheless, accurate approximations can still be attained through appropriate regularization. Currently, the Truncated Singular Value Decomposition (TSVD) is the most commonly used regularization technique, and several fast algorithms have been developed on this basis \cite{Lyon2011,Matthysen2016FAST,zhao2024fast}. However, a well-known drawback of these conventional approaches is the pronounced oscillatory behavior of the solution in the extended region. Such oscillations not only undermine the stability of derivative approximations--critical in solving partial differential equations (PDEs)--but also limit the utility of Fourier extension in broader applications, including signal processing, image reconstruction, and numerical differentiation, where smoothness or physical plausibility outside the original domain is desirable. Indeed, spurious oscillations can introduce artificial high-frequency components that degrade the quality of reconstructions and hinder the interpretability of results in inverse problems and computational physics. A post-processing technique discussed in \cite{2012Sobolev} indicates that extension functions with favorable smoothness properties can be constructed without sacrificing accuracy, highlighting the need for a more intrinsic approach to controlling extension behavior.

According to  the extension theorem of  Sobolev spaces, it is possible to identify extension functions with desirable properties. This raises a question: Can we find an effective method to directly calculate such extension functions?
In fact, when $N$ is finite, $\{e^{\text{i}\frac{k\pi}{T}t}\}_{k=-N}^N$ is a  basis  on [-1,1], that is, the solution to the least squares problem is unique. However, when the TSVD is used, since the singular  vectors with singular values $\sigma_k\leq\epsilon$ (where $\epsilon$ is the truncation parameter) are discarded, the problem actually becomes a multi-solution. TSVD only  yields the solution which minimizes the Euclidean norm of the Fourier coefficients.

For smooth functions, the Fourier coefficients exhibit a characteristic decay rate: the smoother the function, the faster the decay. This observation suggests that the selection among admissible solutions should incorporate \emph{a priori} smoothness information. The idea of enforcing specific solution properties through modified regularization targets has a long history in inverse problems. In his seminal work \cite{tikhonov1963regularization}, Tikhonov (1963) defined the regularization solution $x_{\alpha}$ of the equation
\begin{equation*}
  \mathcal{K}x=y
\end{equation*}
as the solution to the following least-squares problem:
\begin{equation*}
  \min \{\|\mathcal{K}x-y\|^2+\alpha\|\mathcal{L}x\|^2\},
\end{equation*}
where  $\mathcal{L}$ is a differential operator. Along similar lines, Hansen et. al developed GTSVD and MTSVD methods \cite{hansen1989regularization,hansen1992modified,hansen1994regularization} and Natterer (1984) introduced ``regularization in Hilbert scales'' \cite{natterer1984error}.
All these approaches share a common philosophy: the approximation target is altered to reflect additional structural information about the desired solution.

In this paper, we will revisit Fourier extension from  solving compact operator equations with perturbed data. From this perspective, we can more effectively address some concerns previously highlighted in the literature. Instead of approximating the minimum-norm least-squares solution, we propose to approximate a \emph{weighted best-approximation solution}, where the weighting is designed to encode expected smoothness or decay properties of the Fourier coefficients. This leads naturally to a Fourier extension method based on a \emph{weighted generalized inverse}. By choosing appropriate weight operators, we can significantly suppress unwanted oscillations in the extended region while preserving high approximation accuracy on the original interval.
A key advantage of this approach is its ability to yield convergent approximations not only for the function itself but also for its derivatives, which is often crucial in applications such as solving differential equations or computing sensitivities.

\subsection{Overview of the paper}
The remainder of the paper is organized as follows. In Section~2, we formulate the Fourier extension problem as a compact operator equation and introduce the concept of weighted best-approximation solutions. Based on this formulation, we construct a generalized truncated SVD (GTSVD) regularization scheme tailored to Fourier extensions and provide a rigorous convergence analysis under suitable smoothness assumptions. The theoretical results show that, with appropriately chosen weight operators, the proposed method achieves optimal approximation accuracy on the original interval while simultaneously controlling the weighted norm of the coefficient error. This control directly translates into a suppression of spurious oscillations in the extended region, as confirmed by our error estimates. Section~3 discusses the numerical implementation under uniform sampling, including the choice of parameters and the practical design of the weight functions. Numerical experiments presented in Section~4 demonstrate that the proposed method effectively reduces oscillations in the extended domain and yields substantially more accurate approximations of the function and its derivatives. Finally, Section~5 concludes the paper with a summary of the main results and some perspectives for future work.

\section{Formulation of the Fourier extension by solving compact operator equation}
\subsection{A Revisit of the Fourier Extension}
Let $T>1$, $\Lambda_1=[-1,1]$ and $\Lambda_2=(-T,T)$.
Denote by $L^2(\Lambda_i)$ and $H^s(\Lambda_i)$, $(i=1,2)$, the usual Lebesgue
and Sobolev spaces, respectively.
We consider the system obtained by restricting the orthonormal Fourier basis
on $\Lambda_2$ to $\Lambda_1$, namely
\begin{equation*}
  \Phi=\{\phi_{\ell}\}_{\ell\in\mathbb{Z}},\quad
  \phi_0(x)=\frac{1}{\sqrt{2}},\quad
  \phi_{\ell}(x)=e^{\mathrm{i}\frac{\pi}{T}\ell x},\ \ell=\pm1,\pm 2,\ldots,
  \quad x\in \Lambda_1.
\end{equation*}
As shown in \cite{Huybrechs2010}, this system forms a tight frame for
$L^2(\Lambda_1)$.

We define the operators
\begin{equation}
\begin{aligned}
 &\mathcal{F}: \ell^2 \rightarrow L^2(\Lambda_1), \quad
 {c}=\{\hat{c}_{\ell}\}_{\ell=-\infty}^{\infty}
 \mapsto \sum_{\ell=-\infty}^{\infty} \hat{c}_{\ell}\phi_{\ell},\\
 &\mathcal{B}_{h}: \mathcal{D}(\mathcal{B}_{h}) \subset \ell^2 \rightarrow \ell^2,
 \quad {c}\mapsto \{h(|\ell|)\hat{c}_{\ell}\}_{\ell=-\infty}^{\infty},
\end{aligned}
\end{equation}
where $h:[0,\infty)\to(0,\infty)$ is a non-decreasing function satisfying
$\lim_{x\to\infty} h(x)=\infty$.
We then consider the operator equation
\begin{equation}\label{compacteq}
  (\mathcal{F}c)(x)
  :=\sum_{k=-\infty}^{\infty}\hat{c}_k\phi_k(x)
  = f(x),\quad x\in\Lambda_1.
\end{equation}

If $f\in H^p(\Lambda_1)$, then by the Sobolev extension theorem,
equation \eqref{compacteq} admits multiple solutions.
In this setting, the least-squares solution and the best-approximation solution
are defined as follows.
\begin{itemize}
  \item A vector $c^{\star}\in\ell^2$ is called a \emph{least-squares solution}
  of \eqref{compacteq} if
  \begin{equation*}
    \|\mathcal{F}c^{\star}-f\|_{L^2(\Lambda_1)}
    =\inf\{\|\mathcal{F}c-f\|_{L^2(\Lambda_1)} \mid c\in\ell^2\}.
  \end{equation*}
  \item A vector $c^{\dag}\in\ell^2$ is called the
  \emph{best-approximation solution} of \eqref{compacteq} if
  \begin{equation*}
    \|c^{\dag}\|_{\ell^2}
    =\inf\{\|c^{\star}\|_{\ell^2}
    \mid c^{\star}\ \text{is a least-squares solution of \eqref{compacteq}}\}.
  \end{equation*}
\end{itemize}

Due to limitations of machine precision, in practice one can only solve the
perturbed equation corresponding to \eqref{compacteq},
\begin{equation}\label{compacteqnoisy}
  \mathcal{F}c = f^{\epsilon}.
\end{equation}
It is generally assumed that
\begin{equation}
  |f^{\epsilon}(x)-f(x)| \leq \epsilon,
  \quad \forall x\in\Lambda_1,
\end{equation}
which implies
\begin{equation}
  \|f^{\epsilon}-f\|_{L^2(\Lambda_1)} \leq \sqrt{2}\,\epsilon.
\end{equation}
The classical Fourier extension method can be interpreted as computing an
approximation of $c^{\dag}$ via the truncated singular value decomposition
(TSVD).

If the function $f$ possesses a certain degree of smoothness, then
equation \eqref{compacteq} admits solutions whose Fourier coefficients exhibit
corresponding decay properties.
To characterize this behavior, we introduce an index function $h_s$.

\medskip
\noindent\textbf{Assumption.}
There exist a constant $E>0$ and a function $h_s$ such that
\begin{equation}\label{smoothcondition}
  \|\mathcal{B}_{h_s}c\|_{\ell^2} \leq E.
\end{equation}

Under this assumption, it is natural to select the following
weighted best-approximation solution as the approximation target:
\begin{equation*}
  \|c_{\mathcal{B}}^{\dag}\|_{\ell^2}
  =\inf\{\|\mathcal{B}_{h_s}c^{\star}\|_{\ell^2}
  \mid c^{\star}\ \text{is a least-squares solution of \eqref{compacteq}}\}.
\end{equation*}
In practical computations, however, the precise form of $h_s$ is generally
unknown.
Therefore, when designing regularization methods, one typically employs a
prescribed function $h_r$ and analyzes the convergence properties under suitable
relations between $h_s$ and $h_r$.

Let
\begin{equation}
  \mathcal{R}=\mathcal{B}_{h_r}, \qquad
  \mathcal{G}=\mathcal{F}\mathcal{R}^{-1}.
\end{equation}
It follows that $\mathcal{G}$ is a compact operator.
We denote by $(\sigma_i, v_i, u_i)$ a singular system of $\mathcal{G}$, where
\begin{equation}
  \sigma_1 \geq \sigma_2 \geq \cdots \geq \sigma_k \geq \cdots > 0.
\end{equation}

For any $g\in L^2(\Lambda_1)$, we define the truncated inverse operator
$\mathcal{T}_k : L^2(\Lambda_1)\rightarrow \ell^2$ by
\begin{equation}
  \mathcal{T}_k g
  = \sum_{i=1}^{k} \frac{1}{\sigma_i}\,\langle g, u_i \rangle\, v_i .
\end{equation}
Based on this construction, the generalized truncated singular value
decomposition (GTSVD) solution of the perturbed equation
\eqref{compacteqnoisy} is given by
\begin{equation}\label{TSVDsol}
  c^{\epsilon,\eta}
  = \mathcal{R}^{-1}\mathcal{T}_{\eta} f^{\epsilon},
\end{equation}
where the truncation parameter $\eta$ is chosen according to the
discrepancy principle
\begin{equation}\label{parameterchoice}
  \eta = \sup\{\,k \mid \sigma_k > \epsilon\,\}.
\end{equation}
\subsection{Convergence Analysis of Regularized Solution}

We now investigate the convergence properties of the regularized solution
defined in \eqref{TSVDsol}. The analysis is based on the spectral properties
of the compact operator $\mathcal{G}$ and the approximation properties
of truncated Fourier coefficients.

\begin{lemma}\label{lem:G-estimates}
The following operator norm bounds hold:
\begin{equation}\label{operatorest}
\begin{aligned}
&\|\mathcal{T}_k\| \leq \frac{1}{\sigma_k}, \qquad
\|\mathcal{G}\mathcal{T}_k\| \leq 1, \qquad
\|(I-\mathcal{G}\mathcal{T}_{k-1})\mathcal{G}\| \leq \sigma_k,\\
&\|I-\mathcal{G}\mathcal{T}_k\| \leq 1, \qquad
\|I-\mathcal{T}_k\mathcal{G}\| \leq 1.
\end{aligned}
\end{equation}
\end{lemma}

\medskip
To quantify the approximation error induced by truncating the Fourier
coefficients, we introduce a projection operator.

Let $\tau$ be a positive integer to be specified later, and define
\begin{equation}
  \mathcal{P}_{\tau}:\ell^2\rightarrow\ell^2, \qquad
  c \mapsto (\ldots,0,\hat c_{-\tau},\ldots,\hat c_{\tau},0,\ldots).
\end{equation}

\begin{lemma}\label{lem:projection}
Let
\begin{equation}\label{defmid}
  c_{\mathcal{B},\tau} = \mathcal{P}_{\tau} c_{\mathcal{B}}^{\dag},
\end{equation}
then the following estimates hold:
\begin{equation}\label{projest}
\|c_{\mathcal{B}}^{\dag}-c_{\mathcal{B},\tau}\|_{\ell^2}
\leq \frac{1}{h_s(\tau)}E,
\qquad
\|\mathcal{R}c_{\mathcal{B},\tau}\|_{\ell^2}
\leq
\left(\max_{|\ell|\leq\tau}\frac{h_r(|\ell|)}{h_s(|\ell|)}\right)E.
\end{equation}
\end{lemma}

\begin{proof}
Using the monotonicity of $h_s$, we obtain
\[
\|c_{\mathcal{B}}^{\dag}-c_{\mathcal{B},\tau}\|_{\ell^2}^2
=\sum_{|\ell|>\tau}\hat c_\ell^2
\leq \frac{1}{h_s^2(\tau)}
\sum_{|\ell|>\tau}h_s^2(|\ell|)\hat c_\ell^2
\leq \frac{1}{h_s^2(\tau)}\|\mathcal{B}_{h_s}c\|_{\ell^2}^2.
\]
Similarly,
\[
\begin{aligned}
\|\mathcal{R}c_{\mathcal{B},\tau}\|_{\ell^2}^2
&=\sum_{|\ell|\leq\tau}h_r^2(|\ell|)\hat c_\ell^2
=\sum_{|\ell|\leq\tau}
\frac{h_r^2(|\ell|)}{h_s^2(|\ell|)}\,h_s^2(|\ell|)\hat c_\ell^2\\
&\leq
\left(\max_{|\ell|\leq\tau}\frac{h_r^2(|\ell|)}{h_s^2(|\ell|)}\right)
\sum_{|\ell|\leq\tau}h_s^2(|\ell|)\hat c_\ell^2
\leq
\left(\max_{|\ell|\leq\tau}\frac{h_r^2(|\ell|)}{h_s^2(|\ell|)}\right)
\|\mathcal{B}_{h_s}c\|_{\ell^2}^2.
\end{aligned}
\]
\end{proof}

\medskip
Next, we decompose the total error into noise propagation and approximation
errors.

\begin{lemma}\label{lem:mainestimate}
Let
\begin{equation}
  c^{k}= \mathcal{R}^{-1}\mathcal{T}_k f, \qquad
  c^{\epsilon,k}= \mathcal{R}^{-1}\mathcal{T}_k f^{\epsilon}, \qquad
  c^{k}_{\tau}=\mathcal{R}^{-1}\mathcal{T}_k \mathcal{F}c_{\mathcal{B},\tau}.
\end{equation}
Then the following estimates hold:
\begin{equation}
\begin{aligned}
  \|\mathcal{F}\left({c}^{\epsilon, k}-{c}_{\mathcal{B}, \tau}\right)\|\leq \sqrt{2}\epsilon +\frac{1}{h_s(\tau)} E +\sigma_{k+1}\left(\max_{|\ell|\leq \tau}\frac{h_r(|\ell|)}{h_s(|\ell|)}\right) E,\\
   \|\mathcal{R}\left({c}^{\epsilon, k}-{c}_{\mathcal{B}, \tau}\right)\|_{\ell^2}\leq \frac{1}{\sigma_k}\left(\sqrt{2}\epsilon+\frac{1}{h_s(\tau)}E\right)+\left(\max_{|\ell|\leq \tau}\frac{h_r(|\ell|)}{h_s(|\ell|)}\right)  E.
\end{aligned}
\end{equation}
\end{lemma}
\begin{proof}
Using the triangle inequality and lemmas \ref{lem:G-estimates} and \ref{lem:projection}, we have
 \begin{equation}
\begin{aligned}
  \|\mathcal{F}\left({c}^{\epsilon, k}-{c}_{\mathcal{B}, \tau}\right)\|_{L^2(\Lambda_1)}&\leq  \|\mathcal{F}\left({c}^{\epsilon, k}-{c}_{\tau}^{k}\right)\|_{L^2(\Lambda_1)}+\|\mathcal{F}\left({c}_{\mathcal{B},\tau}-{c}_{\tau}^{k}\right)\|_{L^2(\Lambda_1)}\\
  &=\|\mathcal{G}\mathcal{T}_k(f^{\epsilon}-\mathcal{F}{c}_{\mathcal{B},\tau})\|_{L^2(\Lambda_1)}
  +\|(I-\mathcal{G}\mathcal{T}_{k})\mathcal{G}\mathcal{R}{c}_{\mathcal{B},\tau}\|_{L^2(\Lambda_1)}\\
  &\leq \|(f^{\epsilon}-\mathcal{F}{c}_{\mathcal{B},\tau})\|_{L^2(\Lambda_1)}+\sigma_{k+1}\|\mathcal{R}{c}_{\mathcal{B},\tau}\|_{\ell^2}\\
  &\leq \|f^{\epsilon}-f\|_{L^2(\Lambda_1)}+\|f-\mathcal{F}{c}_{\mathcal{B},\tau}\|_{L^2(\Lambda_1)}+\sigma_{k+1}\|\mathcal{R}{c}_{\mathcal{B},\tau}\|_{\ell^2}\\
  &\leq \sqrt{2}\epsilon +\frac{1}{h_s(\tau)} E +\sigma_{k+1}\left(\max_{|\ell|\leq \tau}\frac{h_r(|\ell|)}{h_s(|\ell|)}\right) E,
\end{aligned}
\end{equation}
and
 \begin{equation}
\begin{aligned}
 \|\mathcal{R}\left({c}^{\epsilon, k}-{c}_{\mathcal{B}, \tau}\right)\|_{\ell^2}&\leq
 \|\mathcal{R}\left({c}^{\epsilon, k}-{c}_{\tau}^{k}\right)\|_{\ell^2}+\|\mathcal{R}\left({c}_{\mathcal{B},\tau}-{c}_{\tau}^{k}\right)\|_{\ell^2}\\
 &=\|\mathcal{T}_k(f^{\epsilon}-\mathcal{F}{c}_{\mathcal{B},\tau})\|_{\ell^2}
 +\|(I-\mathcal{T}_{k}\mathcal{G})\mathcal{R}{c}_{\mathcal{B},\tau}\|_{\ell^2}\\
 &\leq \frac{1}{\sigma_k} \|(f^{\epsilon}-\mathcal{F}{c}_{\mathcal{B},\tau})\|_{L^2(\Lambda_1)}+\|\mathcal{R}{c}_{\mathcal{B},\tau}\|_{\ell^2}\\
 &\leq \frac{1}{\sigma_k}\left(\sqrt{2}\epsilon+\frac{1}{h_s(\tau)}E\right)+\left(\max_{|\ell|\leq \tau}\frac{h_r(|\ell|)}{h_s(|\ell|)}\right)  E.
  \end{aligned}
\end{equation}
\end{proof}
\medskip
We are now ready to state the main convergence result.

\begin{theorem}\label{thm1}
Suppose that \eqref{smoothcondition} holds and that
$c^{\epsilon,\eta}$ is defined by \eqref{TSVDsol}-\eqref{parameterchoice}.
If there exists a constant $C_0$ such that
\[
\sup_{\ell}\frac{h_r(|\ell|)}{h_s(|\ell|)}\leq C_0,
\]
then
\begin{equation}\label{errorest1}
\|\mathcal{F}c^{\epsilon,\eta}-f\|_{L^2(\Lambda_1)} = O(\epsilon),
\end{equation}
and
\[
\|\mathcal{R}(c^{\epsilon,\eta}-c_{\mathcal{B}}^{\dag})\|_{\ell^2}=O(1).
\]
\end{theorem}
\begin{proof}
If we choose $\tau$ as
\begin{equation}
  \frac{1}{h_s(\tau)} E=\epsilon,
\end{equation}
then
\begin{equation}
  \begin{aligned}
     \|\mathcal{F}{c}^{\epsilon, \eta}-f\|_{L^2(\Lambda_1)}&\leq \|\mathcal{F}\left({c}^{\epsilon, \eta}-{c}_{\mathcal{B}, \tau}\right)\|_{L^2(\Lambda_1)}+\|\mathcal{F}{c}_{\mathcal{B}, \tau}-f\|_{L^2(\Lambda_1)}\\
     &\leq \sqrt{2}\epsilon+\frac{1}{h_s(\tau)} E +\epsilon\left(\max_{|\ell|\leq \tau}\frac{h_r(|\ell|)}{h_s(|\ell|)}\right) E+\|c_{\mathcal{B}}^{\dag}-{c}_{\mathcal{B}, \tau}\|_{\ell^2}\\
     &\leq \sqrt{2}\epsilon+\frac{1}{h_s(\tau)} E +\epsilon\left(\max_{|\ell|\leq \tau}\frac{h_r(|\ell|)}{h_s(|\ell|)}\right) E+\frac{1}{h_s(\tau)} E\\
     &\leq (\sqrt{2}+2+C_0E)\epsilon,
  \end{aligned}
\end{equation}
and
\begin{equation}
  \begin{aligned}
     \|\mathcal{R}\left({c}^{\epsilon, \eta}-{c}^{\dag}_{\mathcal{B}}\right)\|_{\ell^2}&\leq  \|\mathcal{R}\left({c}^{\epsilon, \eta}-{c}_{\mathcal{B},\tau}\right)\|_{\ell^2}+\|\mathcal{R}\left({c}^{\dag}_{\mathcal{B}}-{c}_{\mathcal{B},\tau}\right)\|_{\ell^2}\\
     &\leq \frac{1}{\sigma_k}\left(\sqrt{2}\epsilon+\frac{1}{h_s(\tau)}E\right)+\left(\max_{|\ell|\leq \tau}\frac{h_r(|\ell|)}{h_s(|\ell|)}\right)  E+\|\mathcal{R}{c}^{\dag}_{\mathcal{B}}\|_{\ell^2}\\
     &\leq  \frac{1}{\epsilon}\left(\sqrt{2}\epsilon+\epsilon\right)+C_0E+C_0E\\
     &\leq \sqrt{2}+1+2C_0E.
  \end{aligned}
\end{equation}
\end{proof}

\begin{corollary}[Convergence rates under typical smoothness assumptions]
Under the assumptions of Theorem~\ref{thm1}, the following convergence
results hold for typical choices of the weight functions $h_s$ and $h_r$.

\begin{enumerate}
\item
If
\[
h_s(x)=\sqrt{1+x^{2s}}, \qquad f\in H^s(\Lambda_1),
\]
and $h_r$ is chosen as
\begin{equation}\label{regu1}
h_r(x)=\sqrt{1+x^{2r}}, \qquad r\leq s,
\end{equation}
then
\[
\|\mathcal{F}c^{\epsilon,\eta}-f\|_{L^2(\Lambda_1)}=O(\epsilon),
\]
and
\[
\|\mathcal{F}c^{\epsilon,\eta}-f\|_{H^r(\Lambda_1)}
\leq
\|\mathcal{R}(c^{\epsilon,\eta}-c_{\mathcal{B}}^{\dag})\|_{\ell^2}
=O(1).
\]
Moreover, by the interpolation inequality in Sobolev spaces, for any
$q\leq r$,
\[
\|\mathcal{F}c^{\epsilon,\eta}-f\|_{H^q(\Lambda_1)}
=O\!\left(\epsilon^{\frac{r-q}{r}}\right).
\]

\item
If
\[
h_s(x)=\mathrm{e}^{c|x|^{s}}, \qquad s\leq 1,
\]
and
\begin{equation}\label{regu2}
h_r(x)=\mathrm{e}^{c|x|^{r}}, \qquad r\leq s,
\end{equation}
then there exists a constant $C_p$ such that, for all $p\geq 0$,
\[
\|\mathcal{F}c^{\epsilon,\eta}-f\|_{L^2(\Lambda_1)}=O(\epsilon),
\]
and
\[
\|\mathcal{F}c^{\epsilon,\eta}-f\|_{H^p(\Lambda_1)}
\leq
C_p\|\mathcal{R}(c^{\epsilon,\eta}-c_{\mathcal{B}}^{\dag})\|_{\ell^2}
=O(1).
\]
Consequently, for any $q\leq p$,
\[
\|\mathcal{F}c^{\epsilon,\eta}-f\|_{H^q(\Lambda_1)}
=O\!\left(\epsilon^{\frac{p-q}{p}}\right).
\]
\end{enumerate}
\end{corollary}

\begin{remark} Some remarks are in order:
\begin{enumerate}
 \item {\bf Relation to the classical Fourier extension method.}
The classical Fourier extension method is recovered as a special case of the
present framework by choosing
\[
h_r(x)\equiv 1.
\]
In this sense, the weighted generalized inverse approach can be viewed as a
natural regularized extension of the classical Fourier extension.

\item {\bf Oversmoothing and parameter choice.}
According to classical regularization theory, oversmoothing does not destroy
convergence provided that the regularization parameter is chosen appropriately.
For example, in \eqref{regu1}, one may formally allow $r>s$. In this case,
to preserve the optimal $L^2$ convergence rate \eqref{errorest1}, the truncation
parameter must satisfy
\[
\sigma_{\eta}\sim \epsilon^{\frac{r}{s}}.
\]
However, such a choice is typically infeasible in practice due to limitations
imposed by machine precision.
Under the practically admissible choice $\sigma_{\eta}\sim\epsilon$, the
oversmoothing case $r>s$ leads instead to the reduced convergence rate
\[
\|\mathcal{F}c^{\epsilon,\eta}-f\|_{L^2(\Lambda_1)}
=O\!\left(\epsilon^{\frac{s}{r}}\right).
\]
An analogous discussion applies to the exponential-weight setting
\eqref{regu2}. This observation highlights the trade-off between stability,
smoothness enforcement, and achievable accuracy in numerical implementations.
\end{enumerate}
\end{remark}

\section{Numerical realization based on uniform sampling}\label{sec3}

\subsection{Discretization based on uniform sampling}
The choice of sampling nodes has a significant impact on the discrete accuracy
of the operator $\mathcal{F}$ and, consequently, on the overall approximation
quality.
Existing studies distinguish between uniform sampling and certain optimal
sampling strategies, such as mapped symmetric Chebyshev nodes
\cite{adcock2014resolution}.
In the present work, we focus on the uniform sampling setting, which is
computationally simple and widely used in practice.

Specifically, the sampling points are chosen as
\begin{equation*}
  t_{j}=\frac{j}{M},\quad j=-M,\ldots,M.
\end{equation*}
For a fixed truncation parameter $N$, we define
\begin{equation*}
  \Phi_N=\{\phi_{\ell}\}_{|\ell|\leq N},
\end{equation*}
and introduce the sampling ratio
\begin{equation*}
  \gamma=\frac{M}{N}.
\end{equation*}
After selecting the parameters $N$, $\gamma$, and $T$, and setting
$L=2T\gamma N$, the discretized counterpart of
\eqref{compacteqnoisy} takes the form
\begin{equation}\label{discreteeq}
  {\bf F}_{\gamma,N}^{T} {\bf c} = {\bf f}^{\epsilon},
\end{equation}
where
\begin{equation}
  \left({\bf F}_{\gamma,N}^{T}\right)_{j,\ell}
  =\frac{1}{\sqrt{L}}\phi_{\ell}(t_{j}),
  \qquad
  \left({\bf f}^{\epsilon}\right)_{j}
  =\frac{1}{\sqrt{L}}f^{\epsilon}(t_{j}).
\end{equation}

In accordance with the continuous formulation
$\mathcal{G}=\mathcal{F}\mathcal{R}^{-1}$, we define the weighted discrete operator
\begin{equation}
  {\bf G}_{\gamma,N}^{T,{\bf R}}
  ={\bf F}_{\gamma,N}^{T}{\bf R}_N,
\end{equation}
where ${\bf R}_N$ is a diagonal matrix with diagonal entries
\begin{equation}\label{diag1}
  ({\bf d})_{\ell}=h_r(|\ell|),\quad \ell=-N,\ldots,N.
\end{equation}

Let
\begin{equation*}
  N_{\gamma}=\min\{2M+1,\,2N+1\}.
\end{equation*}
Then the singular value decomposition of
${\bf G}_{\gamma,N}^{T,{\bf R}}$ is given by
\begin{equation}\label{svdF}
  {\bf G}_{\gamma,N}^{T,{\bf R}}
  =\sum_{i=1}^{N_{\gamma}}{\bf u}_i\,{\bm\sigma}_i\,{\bf v}_i^{T},
\end{equation}
where
\[
  {\bm\sigma}_1\geq \cdots \geq {\bm\sigma}_{N_{\gamma}}\geq 0,
\]
and ${\bf u}_i$, ${\bf v}_i$ denote the left and right singular vectors,
respectively.

Based on this decomposition, an approximation
${\bf c}_N^{\epsilon}$ to the continuous solution
$c^{\epsilon,\eta}$ is defined by
\begin{equation}
  \left({\bf c}_N^{\epsilon}\right)_{\ell}
  =
  \begin{cases}
    \left({\bf R}_N^{-1}\tilde{\bf c}_N^{\epsilon}\right)_{\ell},
      & |\ell|\leq N,\\[2mm]
    0, & |\ell|>N,
  \end{cases}
\end{equation}
where
\begin{equation}\label{GTSVDsoldiscrete}
  \tilde{\bf c}_N^{\epsilon}
  =\sum_{{\bm\sigma}_k>\epsilon}
    \frac{\langle {\bf f}^{\epsilon},{\bf u}_k\rangle}
         {{\bm\sigma}_k}\,{\bf v}_k .
\end{equation}
Finally, the reconstruction
\begin{equation}
  \mathcal{Q}_{\gamma,N}^{T,{\bf R}}{\bf f}^{\epsilon}
  :=\mathcal{F}{\bf c}_N^{\epsilon}
\end{equation}
is used as an approximation of the target function $f$.

\subsection{The design of $h_r$ and the setting of related parameters}
In principle, the weight function $h_r$ may be chosen directly according to
\eqref{regu1} or \eqref{regu2}, reflecting the assumed smoothness of the target
function.
However, for highly oscillatory functions, the corresponding constant $E$ in
the smoothness condition \eqref{smoothcondition} can become excessively large,
which deteriorates the practical approximation accuracy predicted by the
theoretical bounds.

To mitigate this issue, we adopt a modified, $N$-dependent weight function
$h_{r,N}$ that moderates the growth rate of $h_r$ at high frequencies:
\begin{equation}
  h_{r,N}(x)=
  \begin{cases}
    1+\left(\dfrac{x}{C_{r,N}}\right)^{r}, & 0\leq r<\infty,\\[2mm]
    \mathrm{e}^{\frac{x}{C_{r,N}}}-\dfrac{1}{2}, & r=\infty.
  \end{cases}
\end{equation}
Here, the parameter $C_{r,N}$ controls the effective growth rate of the weight
function and thus balances regularization strength and approximation fidelity.
In this work, $C_{r,N}$ is chosen such that
\begin{equation}
  h_{r,N}\!\left(\frac{N}{4}\right)=10,
\end{equation}
which ensures that moderate frequencies are weakly regularized while excessive
penalization of high-frequency components is avoided.

\begin{figure}
			\begin{center}
\subfigure[\label{4a}$\gamma=2$, $p=2$, $N=200$] {
\resizebox*{5.5cm}{!}{\includegraphics{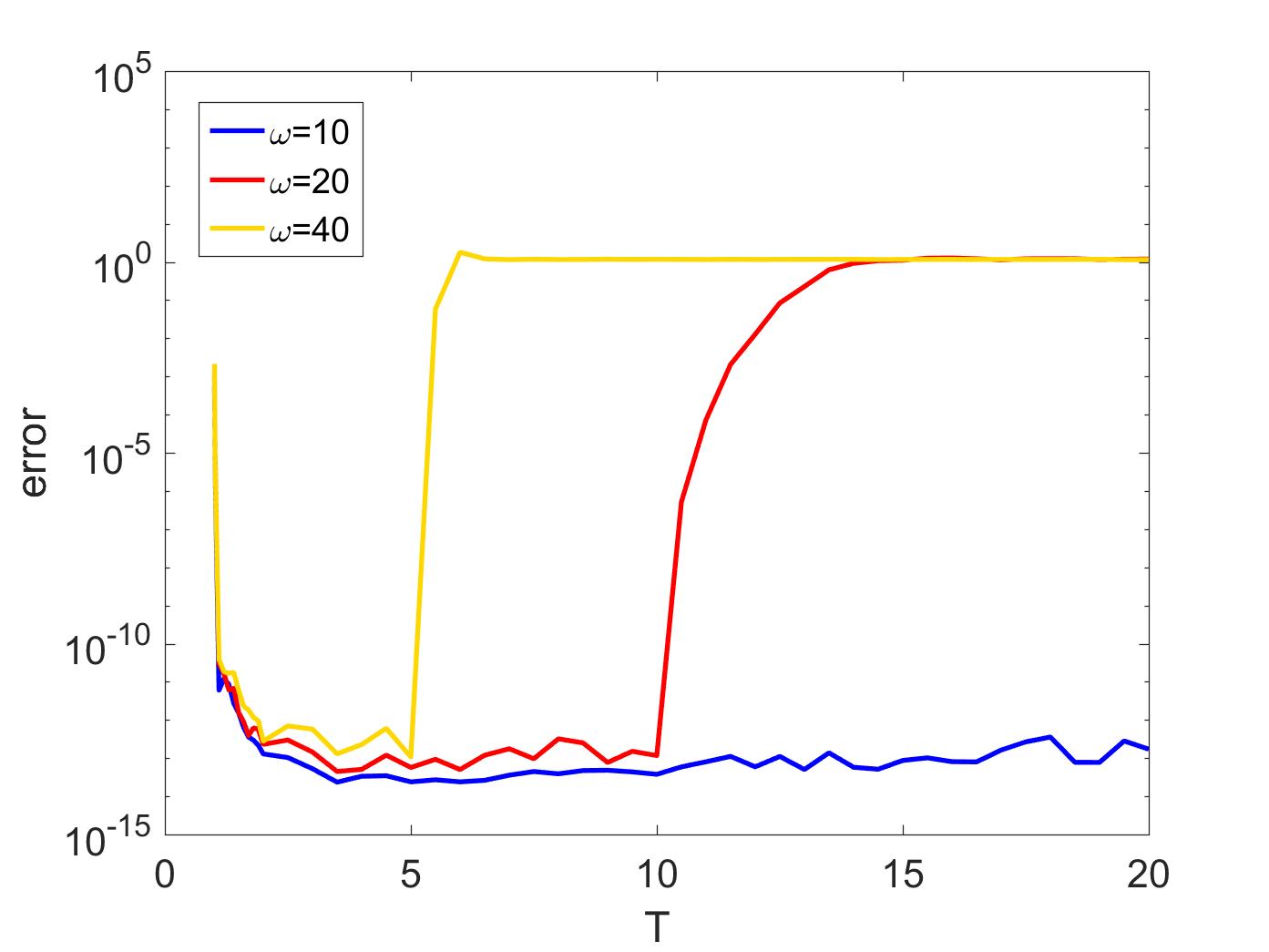}}
}%
\subfigure[\label{4b}$\omega=20$, $p=\infty$,  $N=200$] {
\resizebox*{5.5cm}{!}{\includegraphics{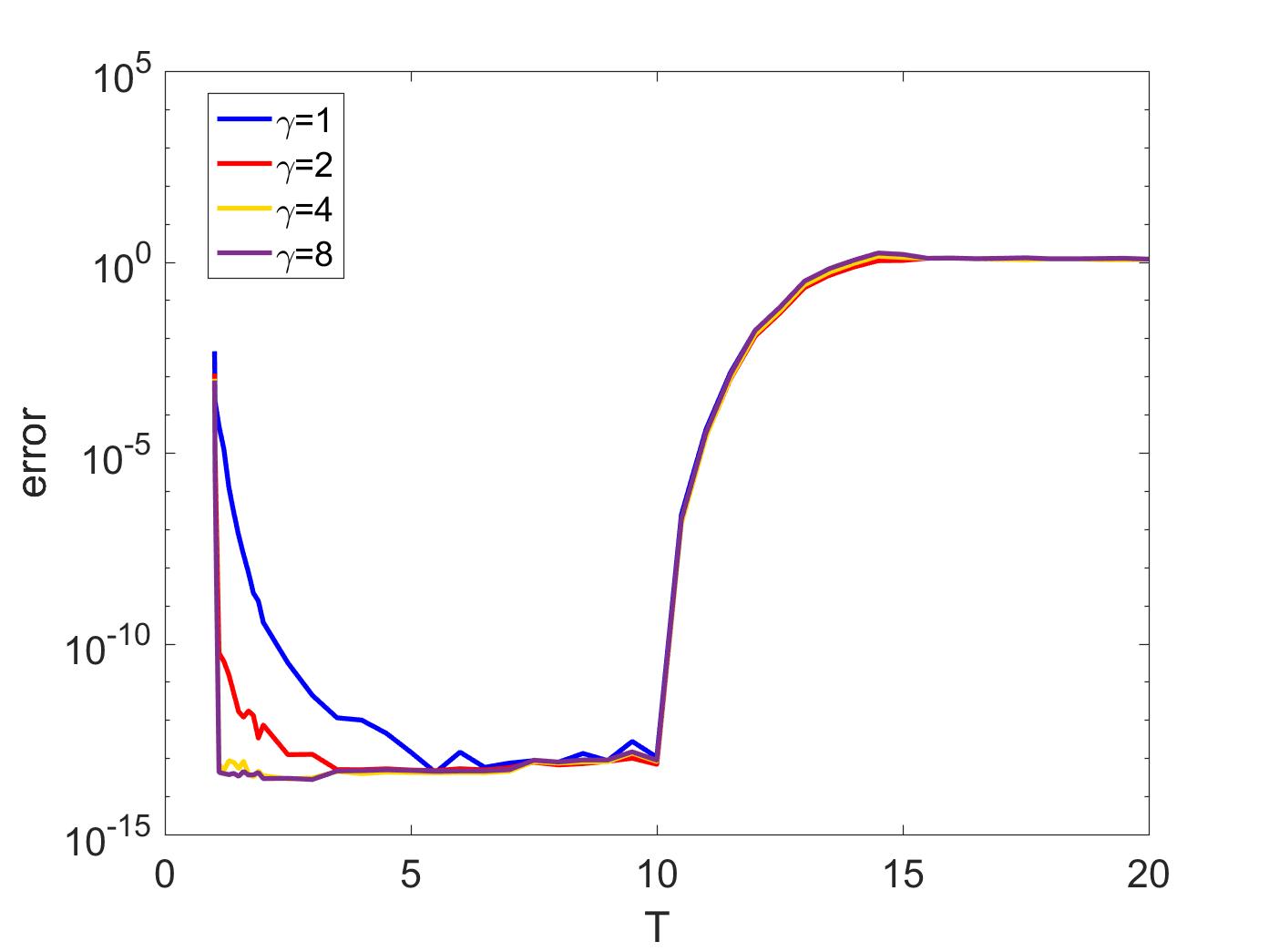}}
}%

 \subfigure[\label{4c}$\omega=20$, $\gamma=2$, $N=200$] {
\resizebox*{5.5cm}{!}{\includegraphics{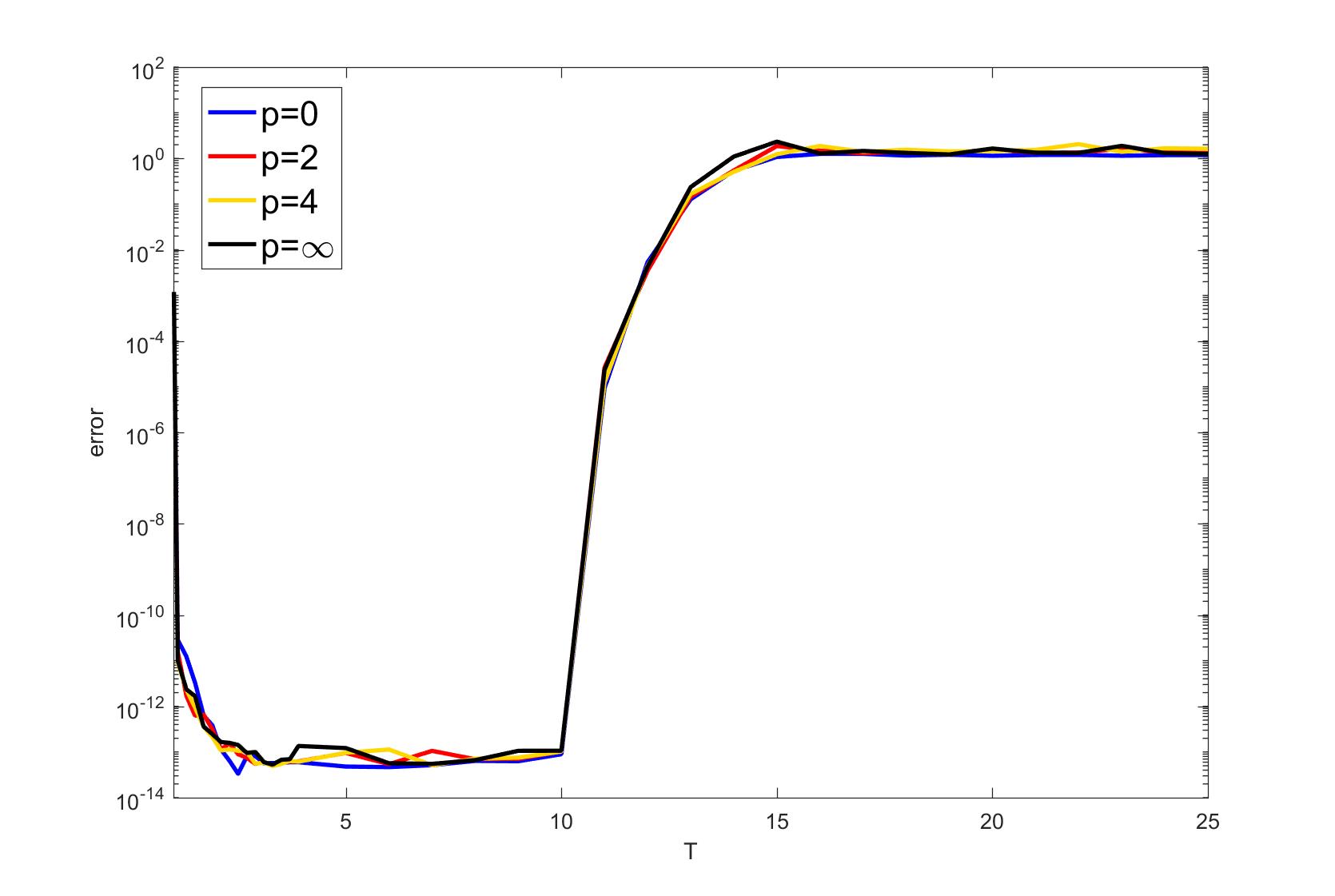}}
}%
 \subfigure[\label{4d}$\omega=20$, $p=2$,$\gamma=2$] {
\resizebox*{5.5cm}{!}{\includegraphics{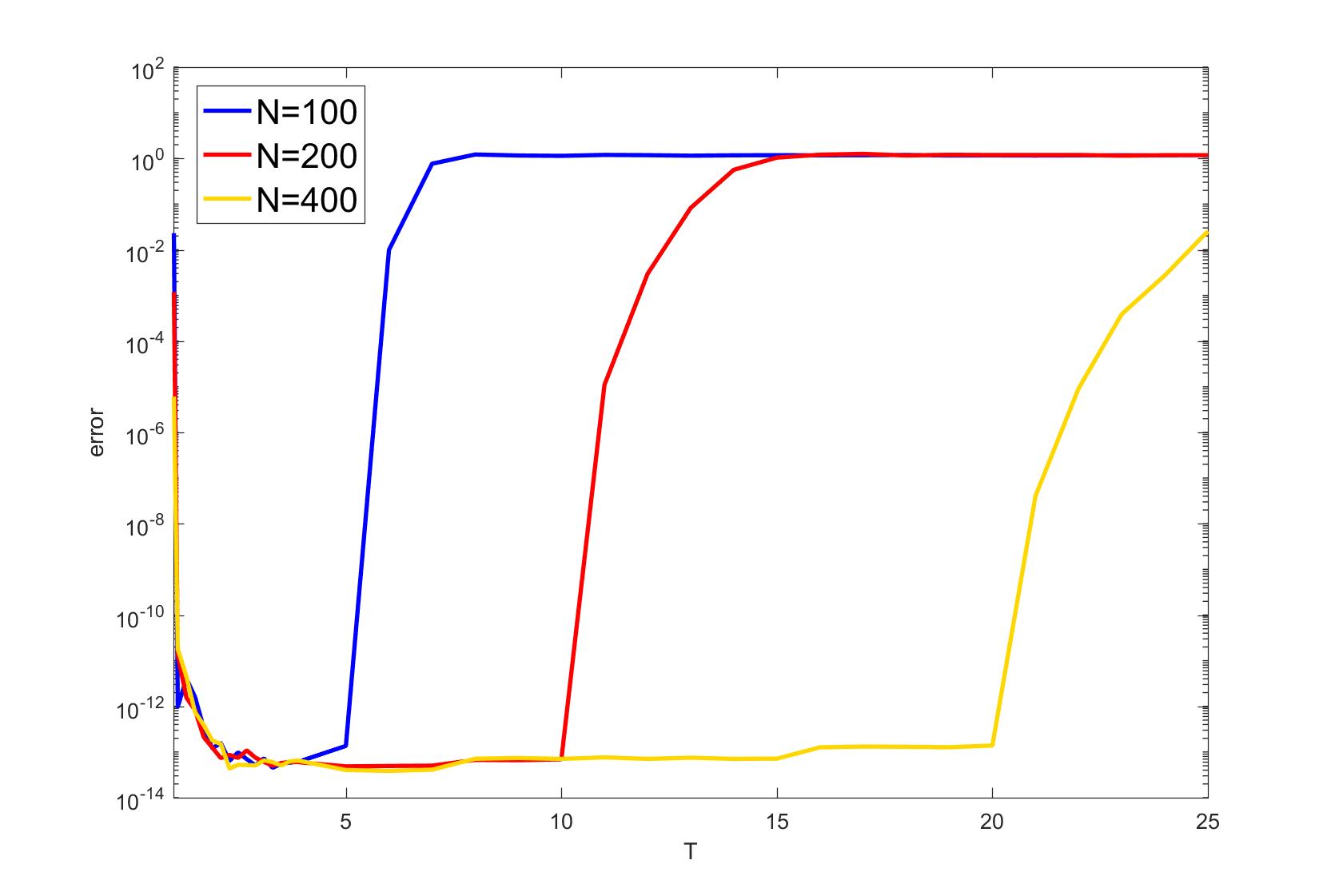}}
}%
		{\caption{The approximation error $\|f-\mathcal{Q}_{\gamma,N}^{T,{\bf R}}{\bf f}^{\epsilon}\|$ against $T$ for different values of $\omega$,  $\gamma$ and $p$. \label{fig4}}}
	\end{center}
\end{figure}

To investigate the influence of the extension parameter $T$, we test the
algorithm using the model function
$f(t)=\exp(\mathrm{i}\pi\omega t)$.
Figure~\ref{fig4} presents the approximation error as a function of $T$ for
different choices of $\omega$, $\gamma$, and $p$.
Three distinct regions can be observed:
\begin{itemize}
  \item Region $I_1:=\{T:1<T<T_1\}$, where the error rapidly decreases to machine
  precision as $T$ increases;
  \item Region $I_2:=\{T:T_1<T<T_2\}$, where the error remains close to machine
  precision;
  \item Region $I_3:=\{T>T_2\}$, where the error increases as $T$ grows.
\end{itemize}

  It can be observed from Fig.~\ref{fig4} that the value of $T_1$ is mainly influenced
by the sampling ratio $\gamma$, while it is essentially independent of the
frequency $\omega$, the truncation parameter $N$, and the regularization order
$p$.
In contrast, the value of $T_2$ depends on $\gamma$, $N$, and $\omega$.

This behavior can be explained as follows.
When all other parameters are fixed, increasing $\gamma$ improves the
discretization accuracy of the matrix ${\bf F}_{\gamma,N}^{T}$ as an
approximation of the continuous operator $\mathcal{F}$.
A similar effect can be achieved by increasing the extension parameter $T$.
Consequently, the approximation accuracy in Region $I_1$ is governed by the
interaction between $\gamma$ and $T$, rather than by either parameter alone.

To further illustrate this observation, we estimate the values of $T_1$ for
different choices of $\gamma$ and report them in Table~\ref{tab1}.
It is seen that the product $T_1\gamma$ remains nearly constant, indicating that
$T_1$ scales approximately like $\gamma^{-1}$.

The value of $T_2$ admits a more direct interpretation.
For the test function $f(t)=\exp(\mathrm{i}\pi\omega t)$, the dominant mode in the
Fourier extension basis $\{\phi_k\}$ corresponds to the index $k=\omega T$.
Therefore, an accurate approximation is possible only when
\[
N \geq \omega T,
\]
which yields the estimate
\begin{equation}\label{T2est}
  T_2 = \frac{N}{\omega}.
\end{equation}

In practical computations, estimate \eqref{T2est} implies that, for functions of
fixed frequency, increasing $T$ requires a corresponding increase in $N$.
Similarly, increasing $\gamma$ leads to a larger number of sampling points $M$.
Both effects increase the computational cost and thus affect the overall
efficiency of the algorithm.
A detailed study of fast implementations is beyond the scope of the present
paper.
For fast algorithms based on boundary interval data, we refer to
\cite{zhao2024fast}, and to the interval partitioning strategies proposed in
\cite{zhao2025local,zhao2025fast}.

  \begin{table}
\begin{center}
  \caption{The approximation value of ${T}_{1}$  for various $\gamma$. \label{tab1} }
  \small
{\begin{tabular*}{\textwidth}{@{\extracolsep\fill}cccccccccccccc} \toprule
&$\gamma=1$&$\gamma=2$&$\gamma=3$&$\gamma=4$&$\gamma=8$\\\midrule
$T_1$&$5.6$&$2.3$&$1.6$&$ 1.2$ &$1.08$\\
\bottomrule
  \end{tabular*}}
\end{center}
\end{table}

\begin{figure}[H]
			\begin{center}
\subfigure[\label{5a}   $f_1(t)$] {
\resizebox*{5.5cm}{!}{\includegraphics{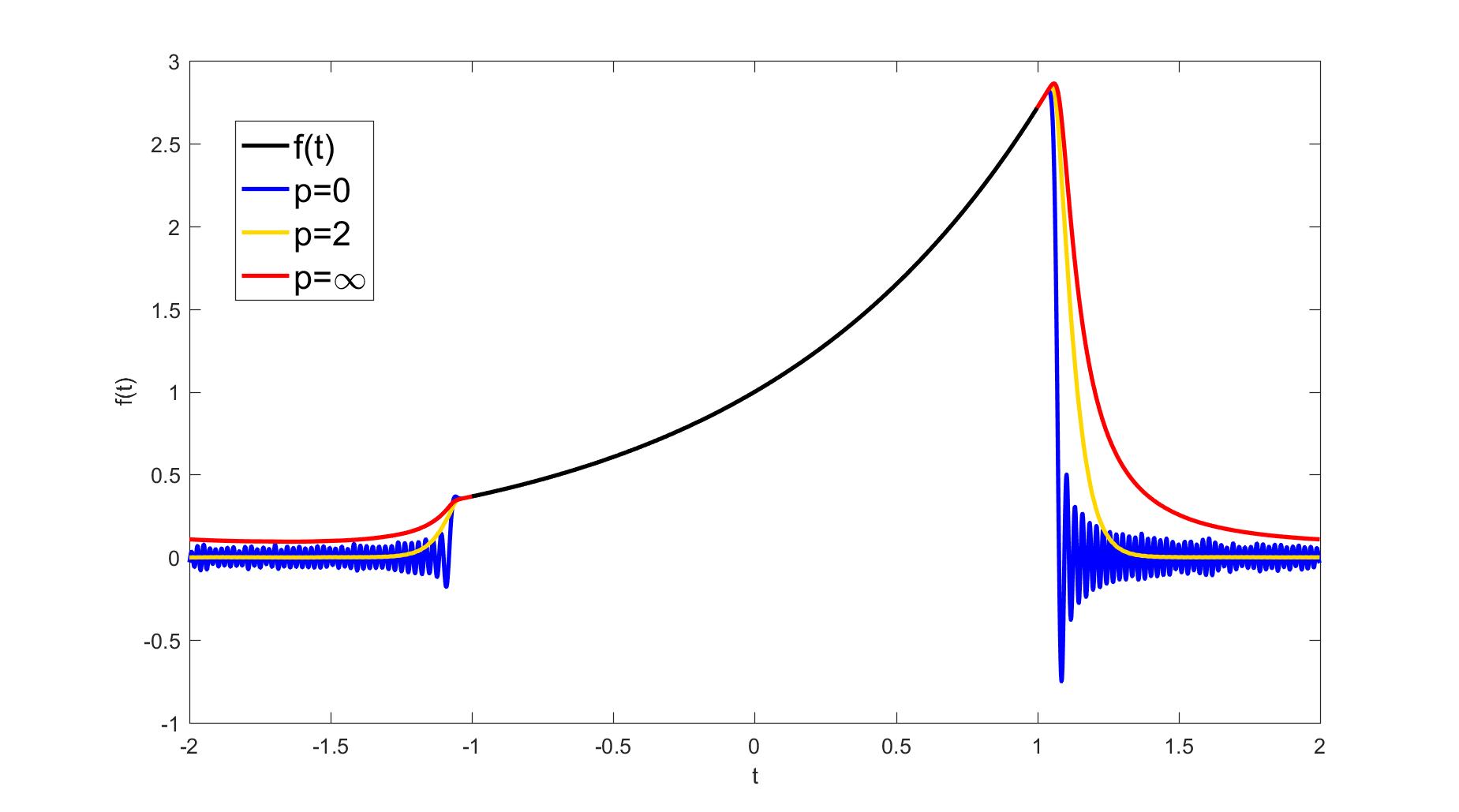}}
}%
\subfigure[\label{5b}  $f_2(t)$] {
\resizebox*{5.5cm}{!}{\includegraphics{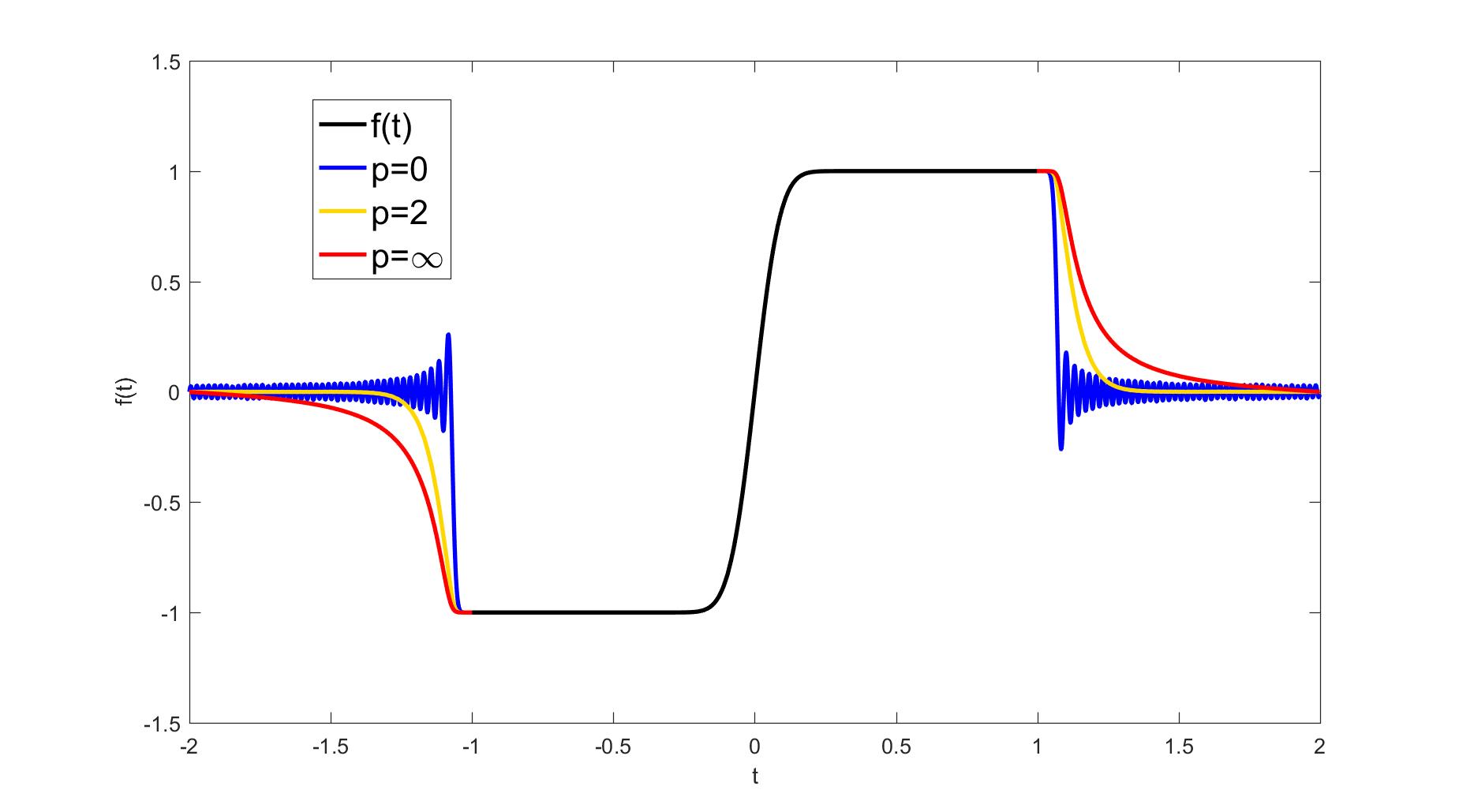}}
}%

\subfigure[\label{5c}  $f_3(t)$] {
\resizebox*{5.5cm}{!}{\includegraphics{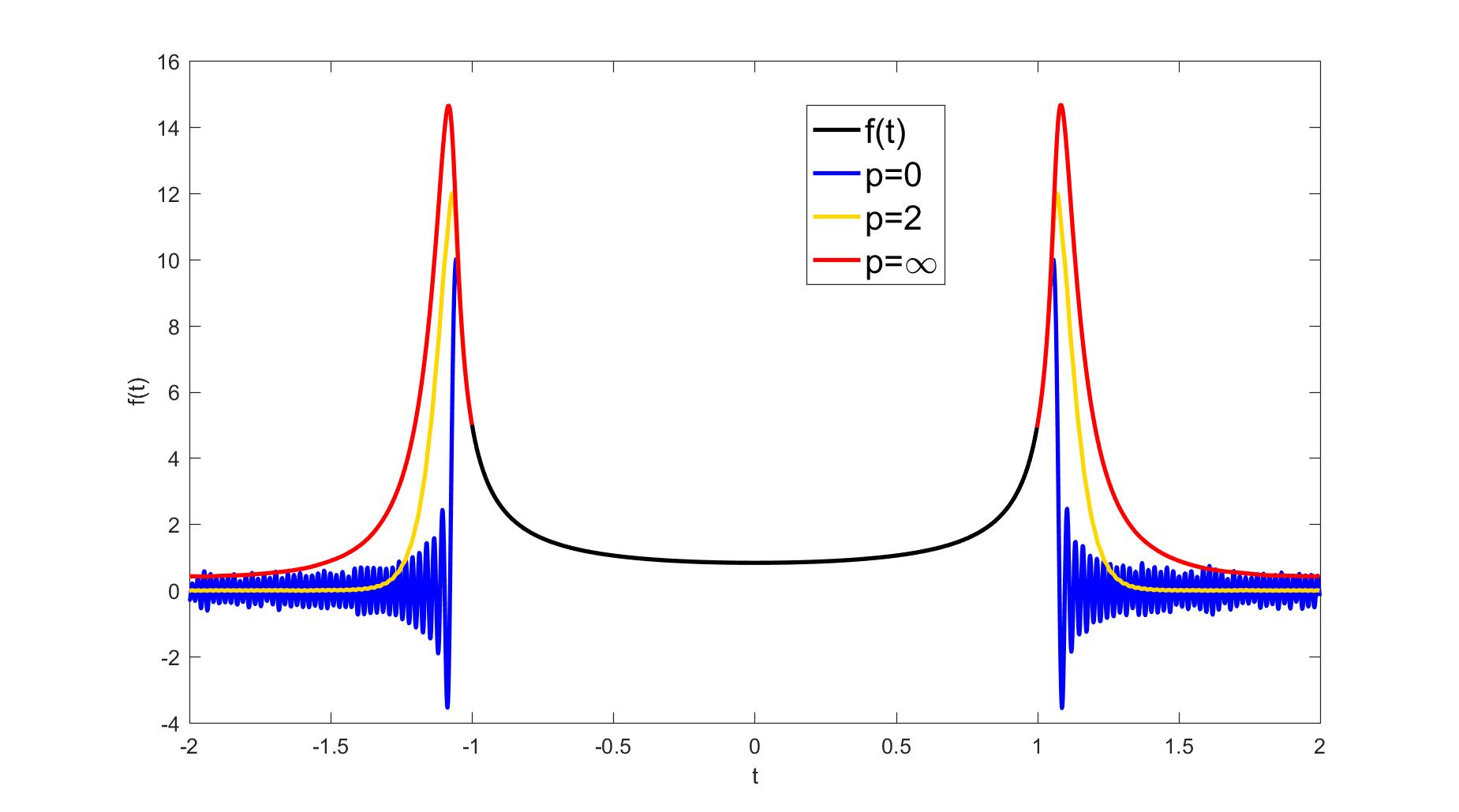}}
}%
\subfigure[\label{5d}   $f_4(t)$] {
\resizebox*{5.5cm}{!}{\includegraphics{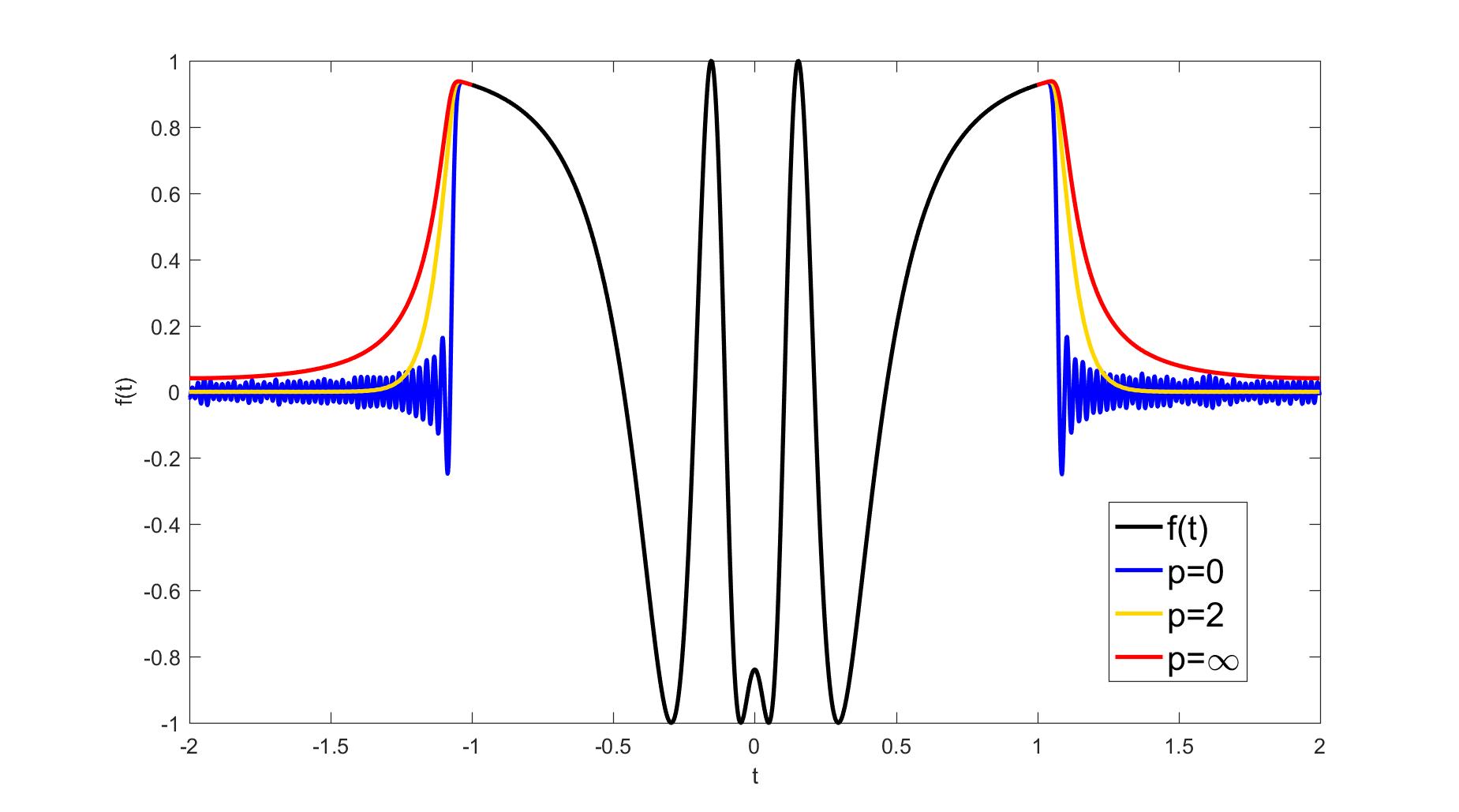}}
}%

	{\caption{ Calculation results of $\mathcal{Q}_{\gamma,N}^{T,{\bf R}}{\bf f}^{\epsilon}$ for the
various test functions ($N=200$).}\label{Fig5}}
	\end{center}
\end{figure}
\section{Numerical tests}\label{sec4}
In this section, we present several numerical experiments to validate the
effectiveness of the proposed method.
All experiments were performed on a Windows~10 system equipped with
16~GB memory and an Intel(R) Core(TM) i7--8500U CPU @~1.80~GHz,
using MATLAB~2016b.

The primary objective of the proposed method is to suppress the spurious
oscillations of Fourier extension approximations in the extended region.
As discussed in the previous section, the parameters $\gamma$ and $T$ have
only a mild influence on this behavior.
Therefore, in all subsequent experiments, we fix
\[
\gamma=3, \qquad T=2,
\]
and focus on the influence of  $p$ and  $N$.

We consider the following four test functions, which exhibit different
smoothness and oscillatory characteristics:
\begin{equation*}
\begin{aligned}
 & f_1(t)=\exp(t), \qquad
   f_2(t)=\mathrm{erf}(10t), \\
 & f_3(t)=\frac{1}{1.2-t^2}, \qquad
   f_4(t)=\cos\!\left(\frac{10}{1+25t^2}\right).
\end{aligned}
\end{equation*}

Figure~\ref{Fig5} displays the weighted projection
$\mathcal{Q}_{\gamma,N}^{T,{\bf R}}{\bf f}^{\epsilon}$ for the test functions
$f_i$, $i=1,2,3,4$, under different choices of the regularization parameter $p$.
According to the theoretical analysis, the case $p=0$ corresponds to the
classical Fourier extension method.
As expected, the corresponding approximations exhibit pronounced oscillations
around zero in the extended region.
For $p>0$, the extension behavior is significantly altered: the oscillations
are strongly suppressed and the extended part of the approximation reflects,
to some extent, the smoothness properties of the target function.
These observations confirm that the proposed weighted regularization
effectively reduces high-frequency oscillations in the extension region
without degrading the approximation quality on the original interval.

\begin{figure}
			\begin{center}
\subfigure[\label{6a} $f_1(t)$] {
\resizebox*{5.5cm}{!}{\includegraphics{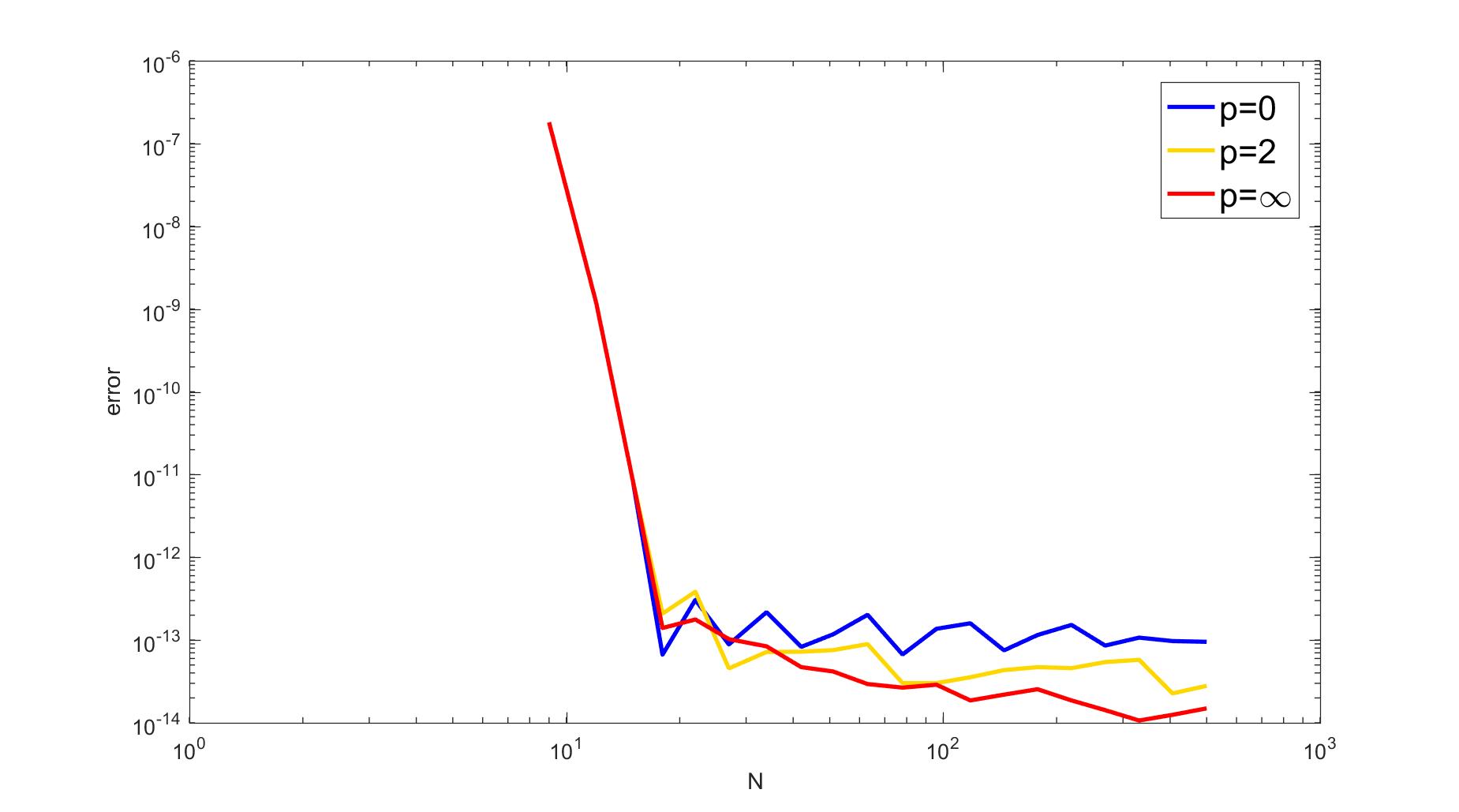}}
}%
\subfigure[\label{6b} $f_2(t)$] {
\resizebox*{5.5cm}{!}{\includegraphics{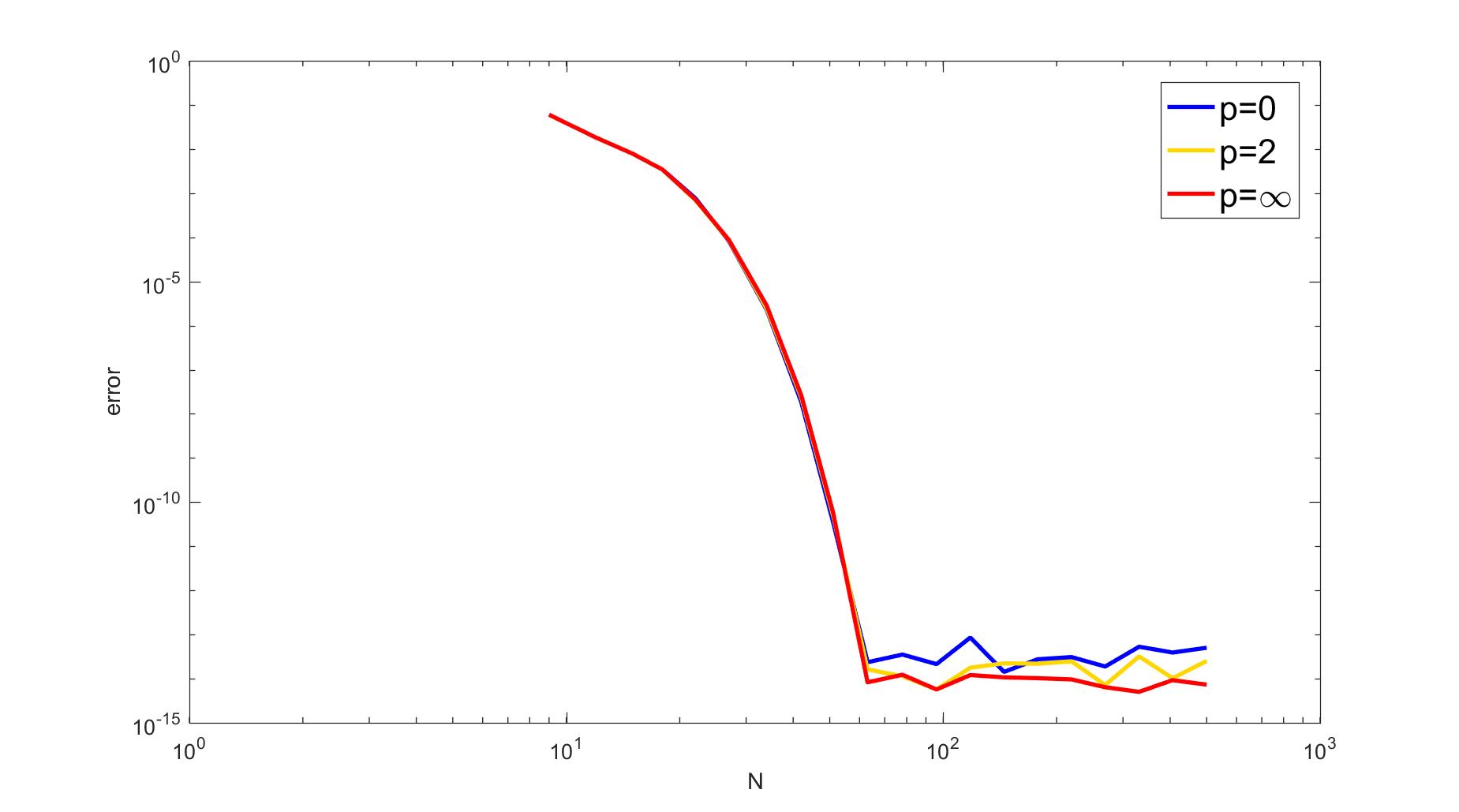}}
}%

\subfigure[\label{6c} $f_3(t)$.] {
\resizebox*{5.5cm}{!}{\includegraphics{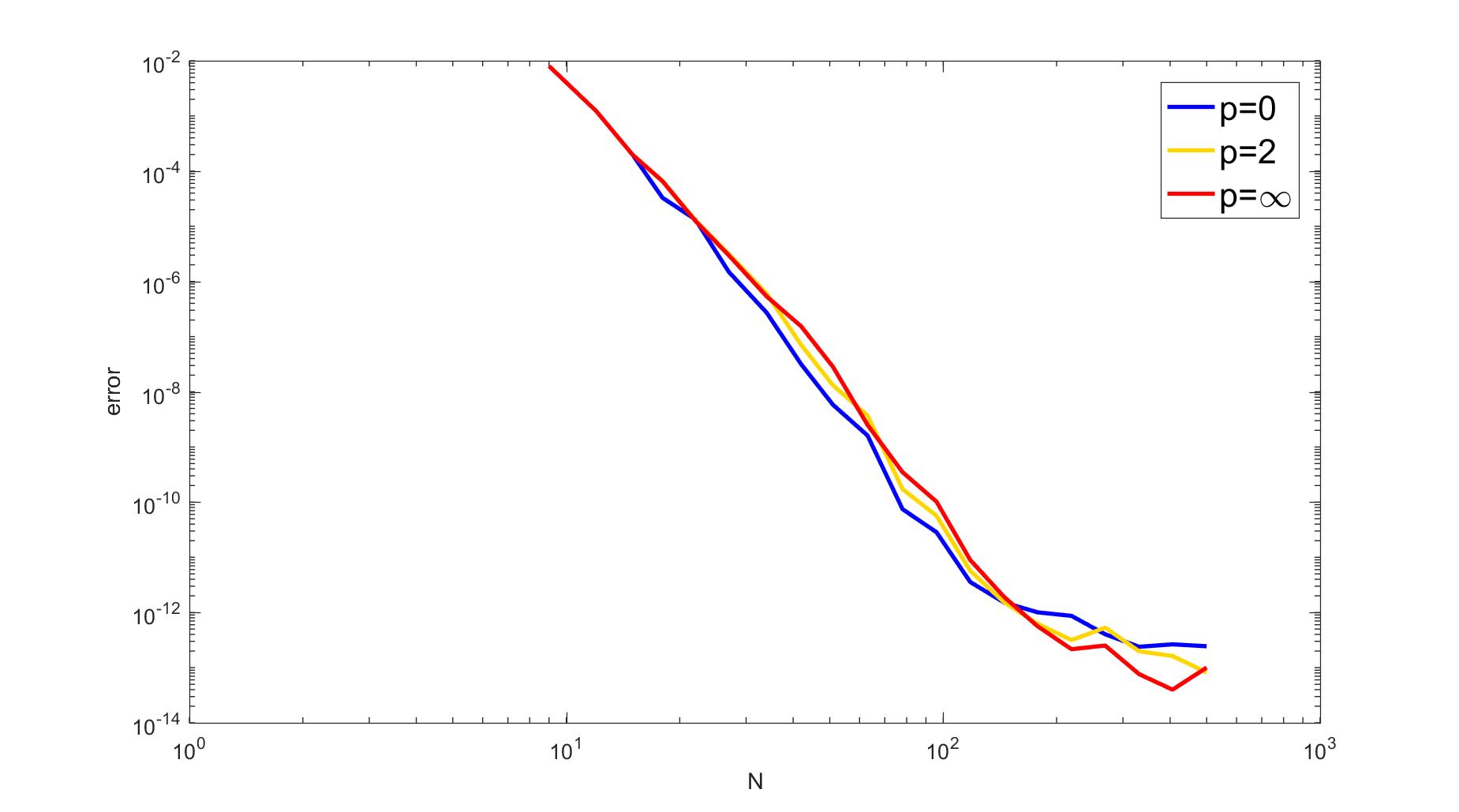}}
}%
\subfigure[\label{6d}$f_4(t)$] {
\resizebox*{5.5cm}{!}{\includegraphics{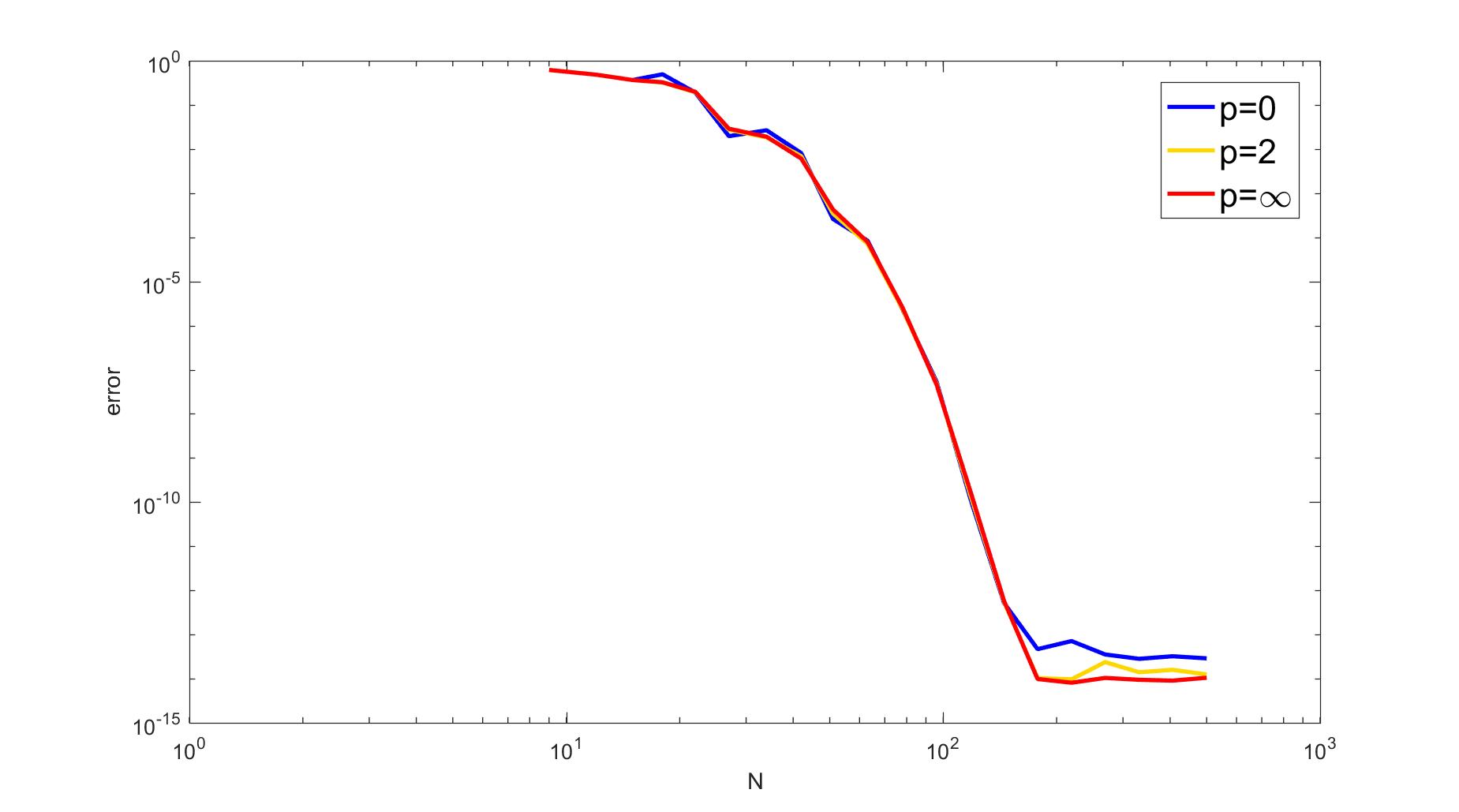}}
}%
	{\caption{ Approximation errors for the various test functions.}\label{Fig6}}
	\end{center}
\end{figure}

\begin{figure}
			\begin{center}
\subfigure[\label{7a} $f_1(t)$] {
\resizebox*{5.5cm}{!}{\includegraphics{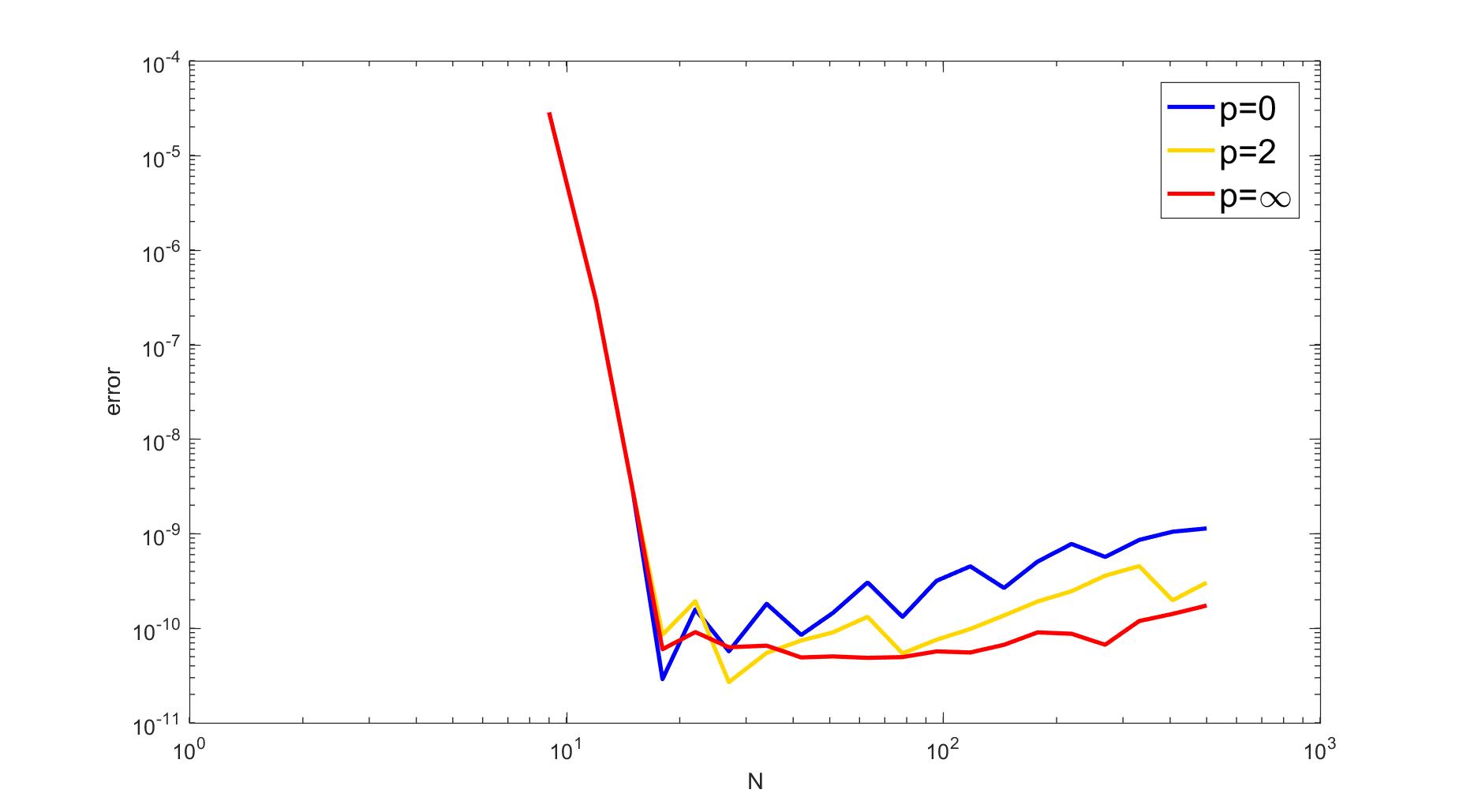}}
}%
\subfigure[\label{7b}$f_2(t)$] {
\resizebox*{5.5cm}{!}{\includegraphics{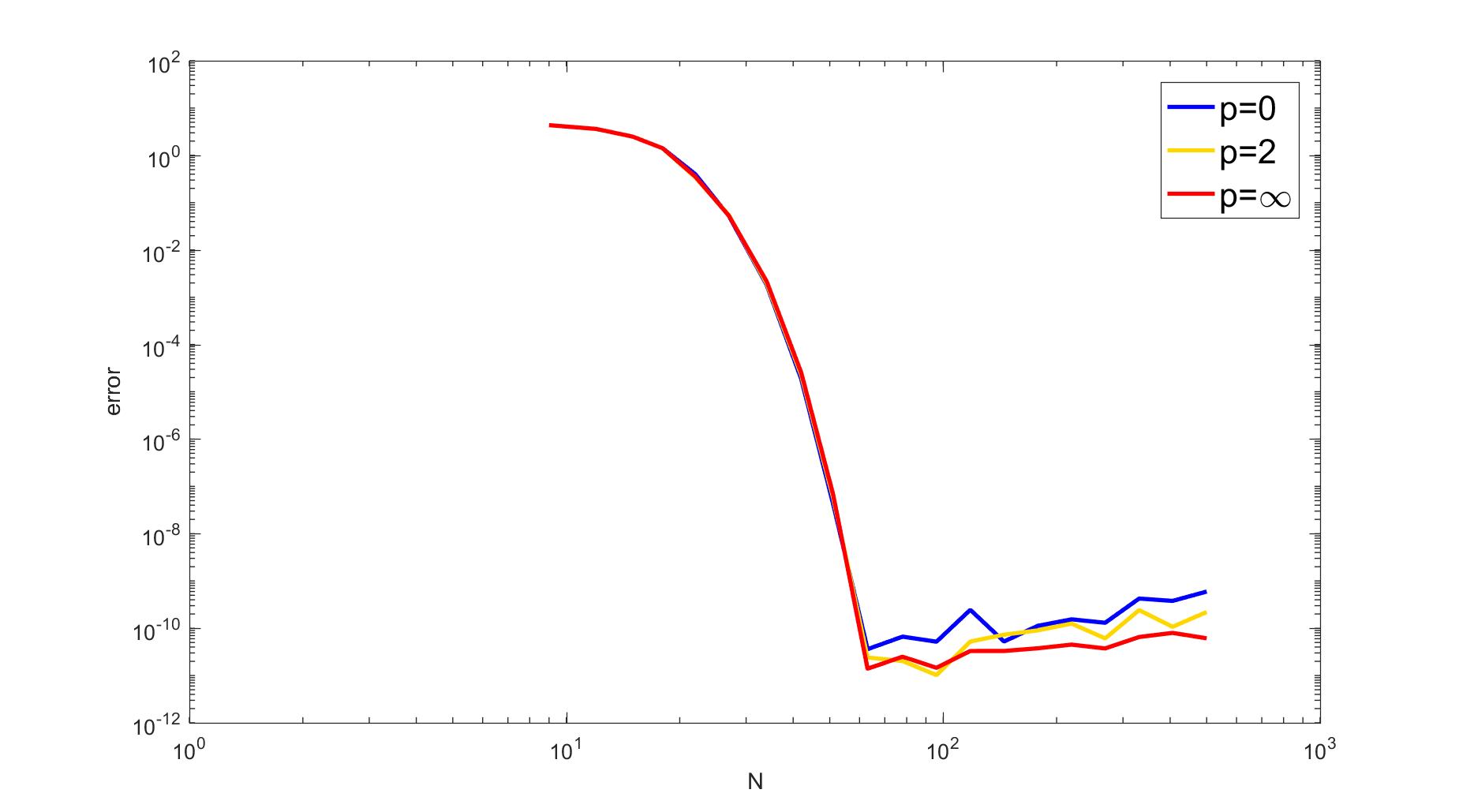}}
}%

\subfigure[\label{7c} $f_3(t)$] {
\resizebox*{5.5cm}{!}{\includegraphics{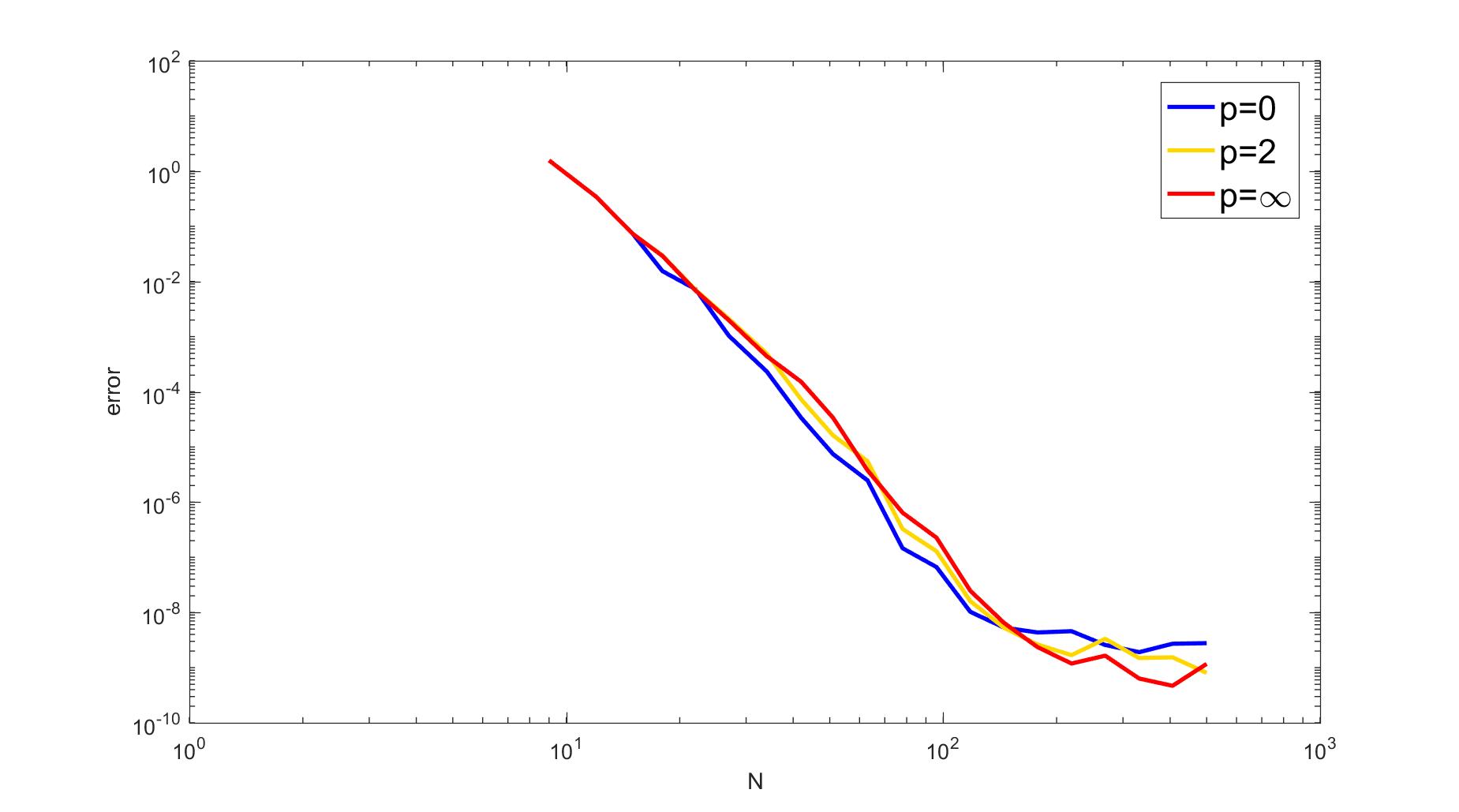}}
}%
\subfigure[\label{7d} $f_4(t)$] {
\resizebox*{5.5cm}{!}{\includegraphics{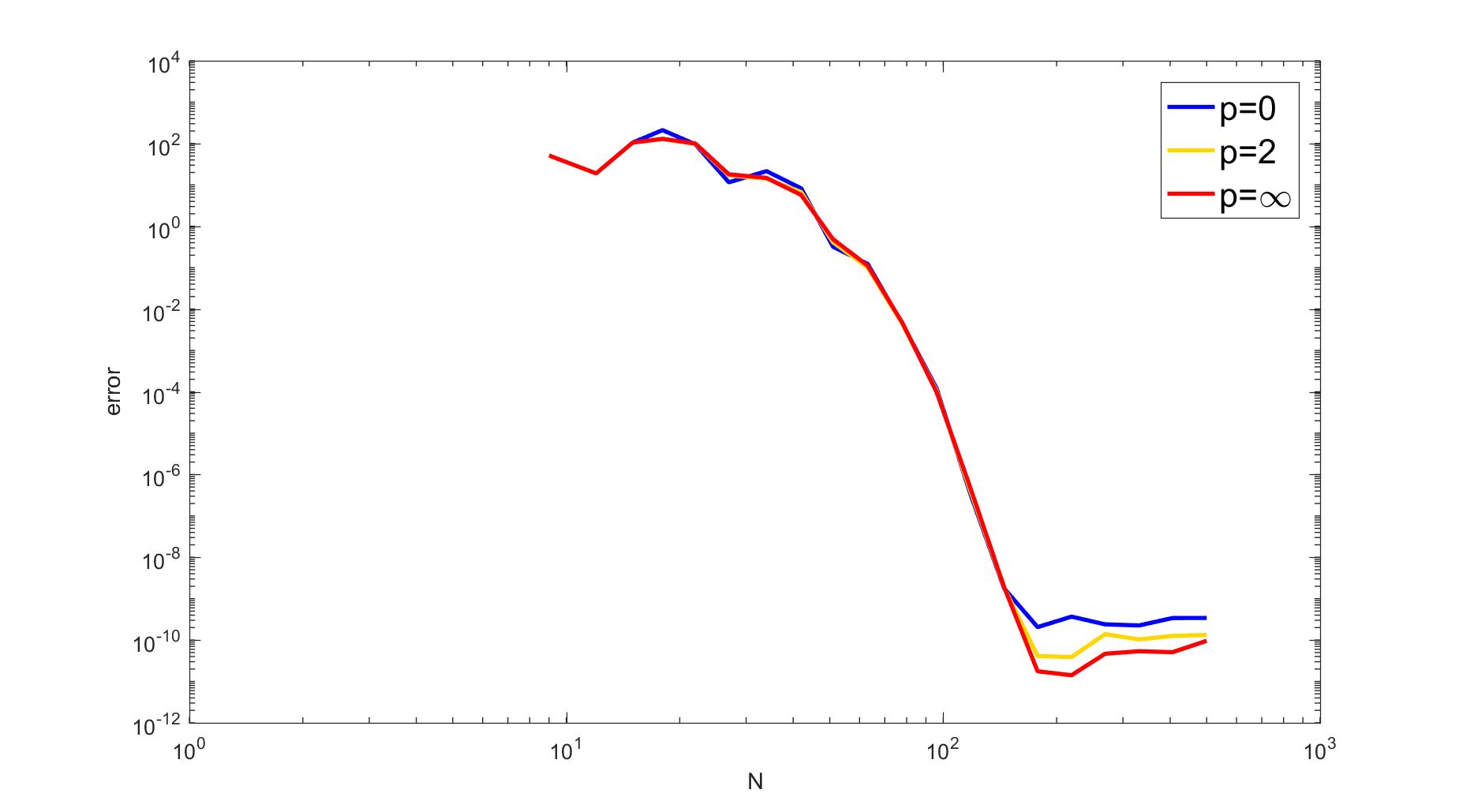}}
}%

	{\caption{Approximation error of the first-order derivative for various test functions.}\label{Fig7}}
	\end{center}
\end{figure}

\begin{figure}
			\begin{center}
\subfigure[\label{8a}$f_1(t)$] {
\resizebox*{5.5cm}{!}{\includegraphics{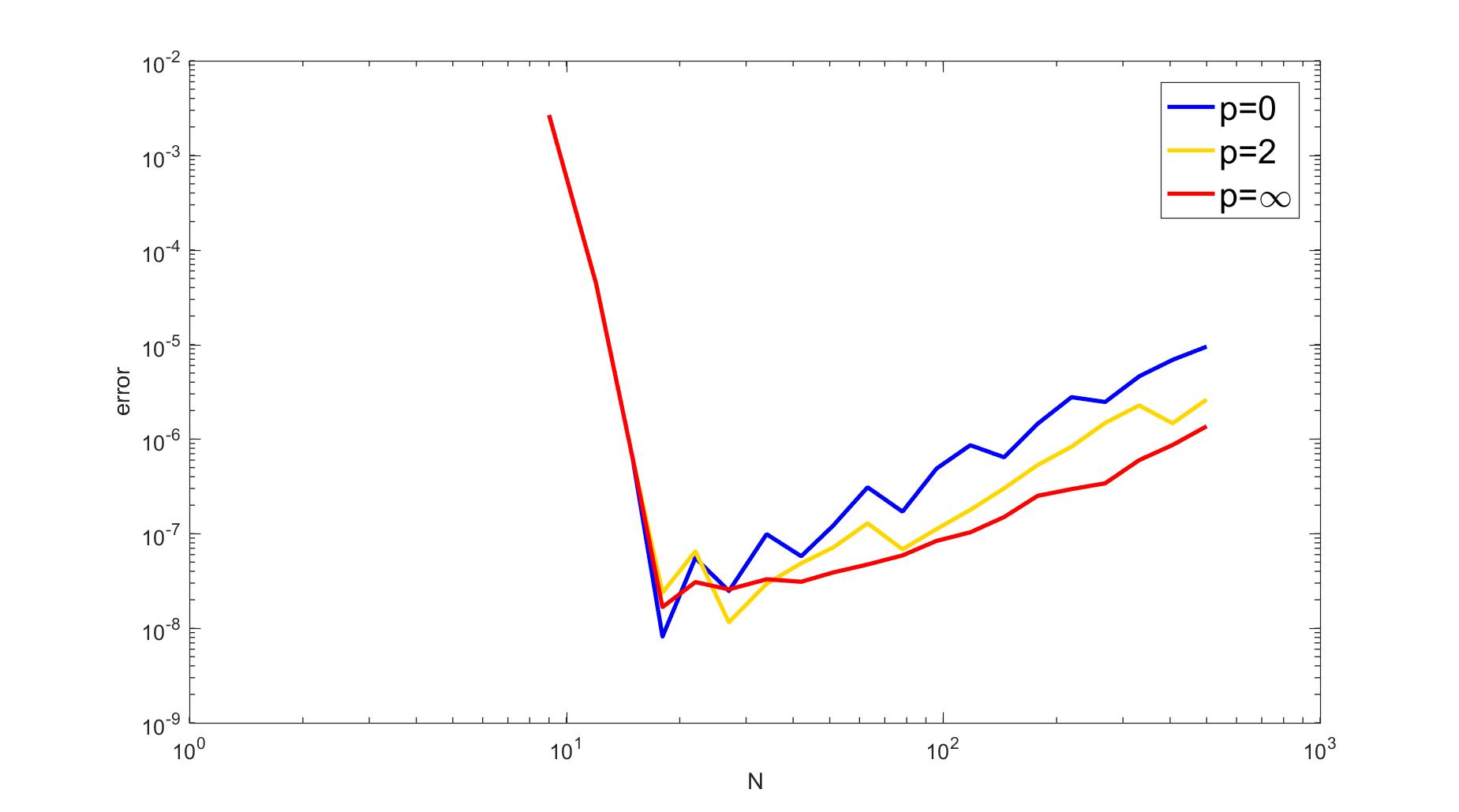}}
}%
\subfigure[\label{8b} $f_2(t)$] {
\resizebox*{5.5cm}{!}{\includegraphics{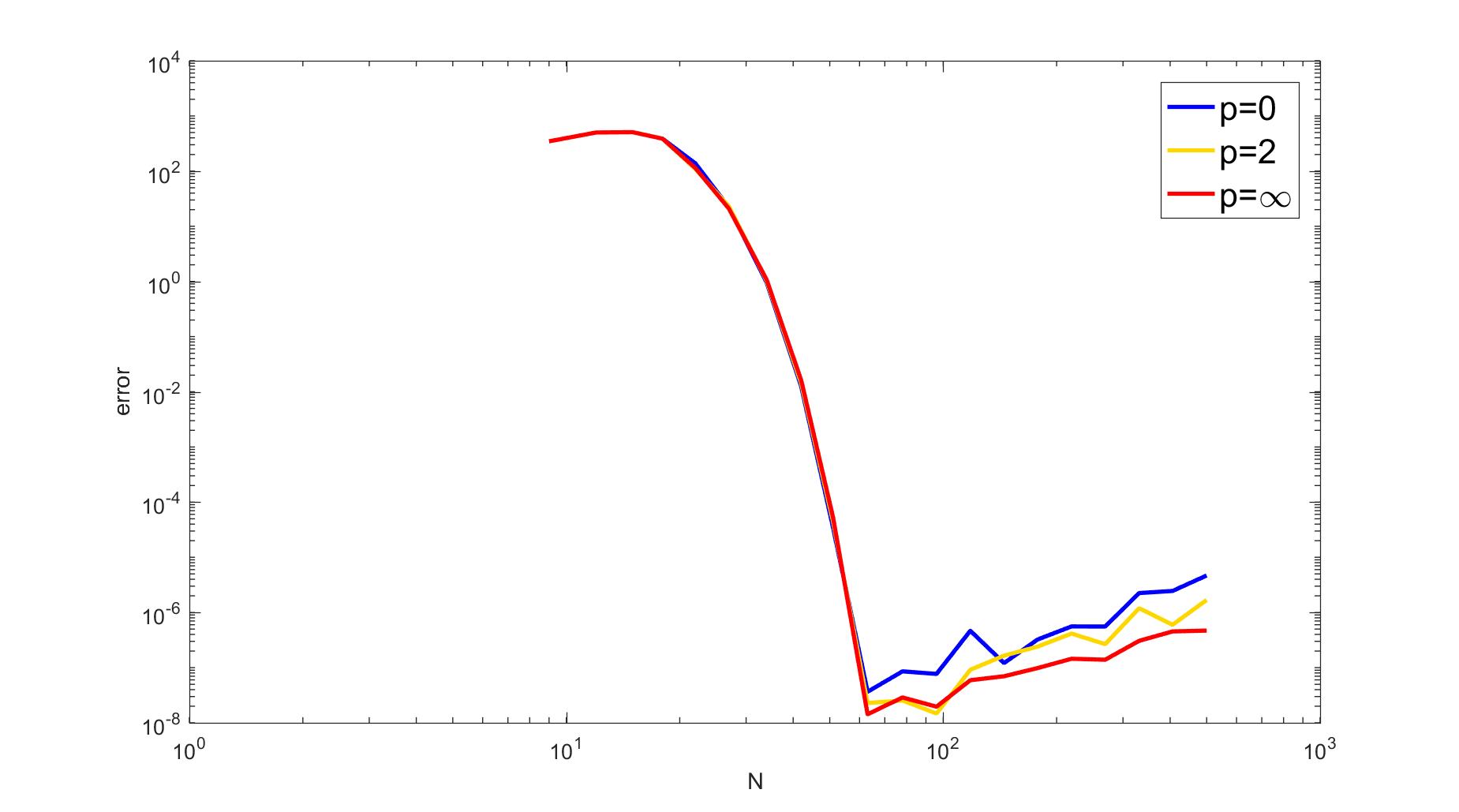}}
}%

\subfigure[\label{8c}$f_3(t)$] {
\resizebox*{5.5cm}{!}{\includegraphics{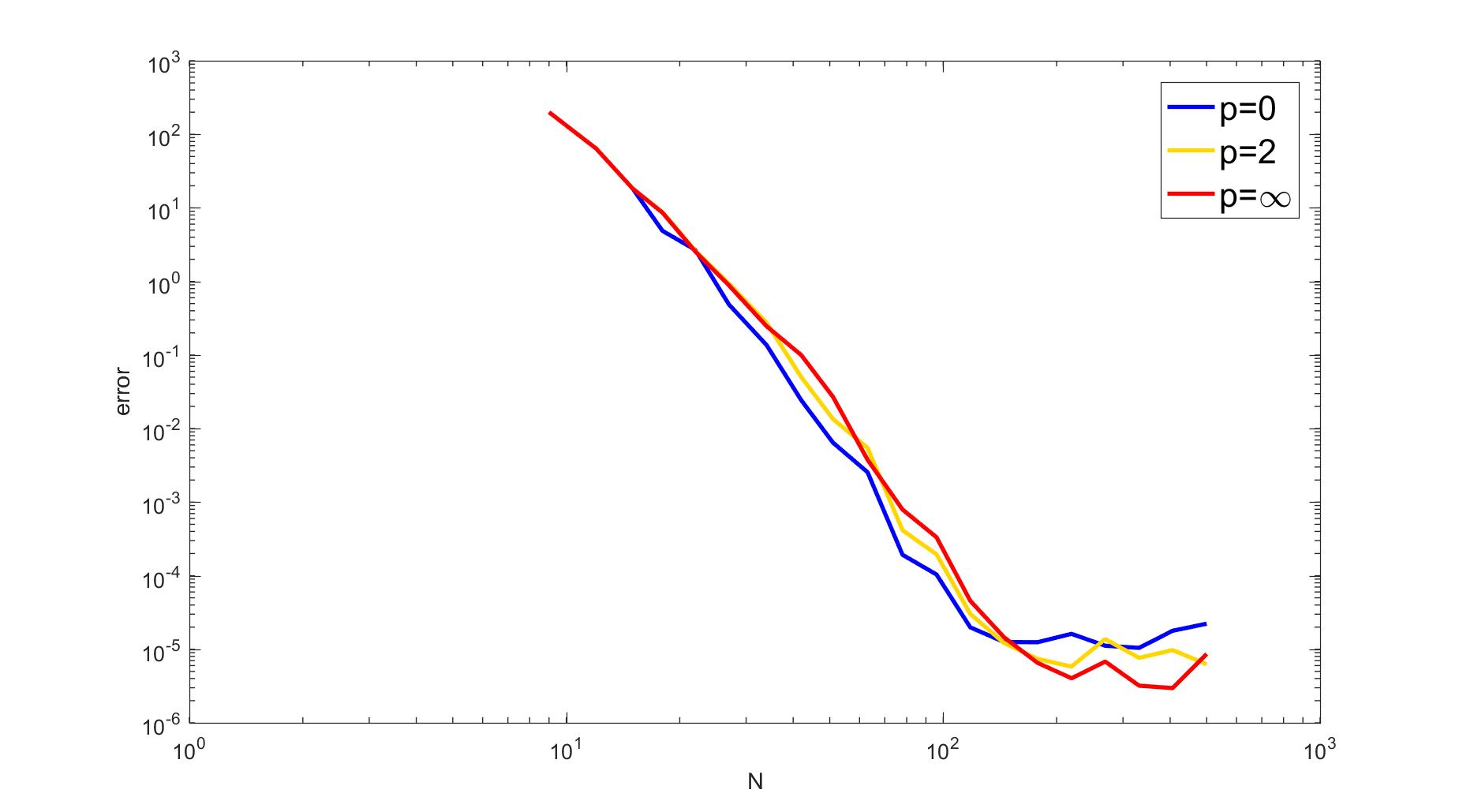}}
}%
\subfigure[\label{8d} $f_4(t)$] {
\resizebox*{5.5cm}{!}{\includegraphics{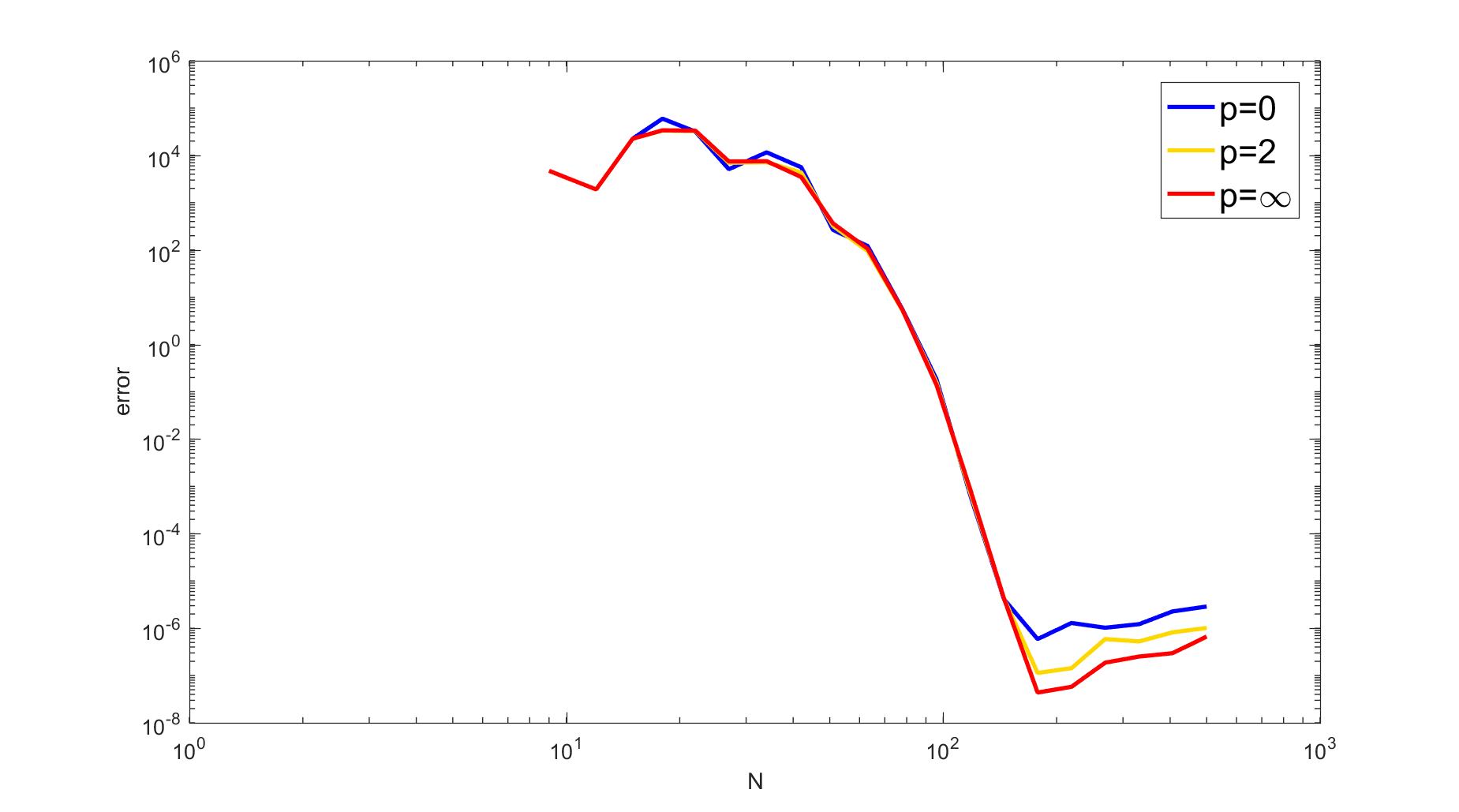}}
}%

	{\caption{Approximation error of the second-order derivative for various test functions.}\label{Fig8}}
	\end{center}
\end{figure}

Figures~\ref{Fig6}--\ref{Fig8} illustrate the dependence of the approximation
errors on the truncation parameter $N$ for the test functions and for their
first- and second-order derivatives.
The derivatives are computed by applying trigonometric interpolation on
$[-T,T]$ using the FFT.
It can be observed that, for the function values themselves, the approximation
errors corresponding to different values of $p$ are nearly indistinguishable.
In contrast, for the first- and second-order derivatives, the proposed method
yields a substantial error reduction compared with the classical Fourier
extension.

We also observe that, once $N$ exceeds a certain threshold, the errors in the
derivative approximations begin to increase as $N$ grows.
This behavior is primarily caused by the derivative evaluation procedure based
on trigonometric interpolation, rather than by the Fourier extension
approximation itself.

\section{Conclusions and remarks}\label{sec5}
In this paper, we proposed a weighted generalized inverse framework for Fourier
extensions and provided a systematic theoretical and numerical investigation of
its properties.
The main conclusions and perspectives can be summarized as follows.

\begin{itemize}
  \item
  From a theoretical point of view, the proposed method constitutes a
  substantial extension of the classical Fourier extension framework.
  By introducing suitable weighted regularization operators and analyzing the
  associated generalized truncated SVD, we established not only convergence of
  the approximation in the $L^2$ sense, but also stability and convergence results
  for derivatives.
  This represents an important improvement over standard Fourier extension
  methods, for which derivative convergence is typically difficult to guarantee
  due to oscillations in the extended region.

  \item
  The analysis in this work demonstrates that smooth extensions can be constructed
  in a controlled manner through appropriate weighting, and that the resulting
  approximations retain high accuracy on the original interval while exhibiting
  improved regularity outside it.
  In this sense, the present framework provides a rigorous justification for
  smoothing strategies that have previously been used mainly from a heuristic or
  numerical perspective.

  \item
  The idea of piecewise approximation, motivated by the sensitivity of global
  Fourier extensions to local singularities, is closely related to the present
  work.
  While existing Fourier extension methods aim at global approximations, local
  irregularities can significantly degrade the overall accuracy.
  This issue has been addressed in our recent work \cite{zhao2025local}, where a
  piecewise Fourier extension strategy was developed and analyzed.
  The weighted regularization framework proposed here provides additional
  theoretical insight into why such localized approaches can be effective, and
  offers a complementary perspective on smooth extension and stabilization.

  \item
  Regarding extensions to higher-dimensional problems, the present paper focuses
  primarily on the theoretical feasibility of constructing smooth and stable
  extensions.
  For practical high-dimensional implementations, we plan to combine the ideas
  developed here with domain decomposition and partitioning strategies, such as
  those proposed in \cite{zhao2025local}, which have demonstrated significant
  advantages in computational efficiency.
  This combination offers a promising pathway for extending the proposed theory
  to multi-dimensional settings.

  \item
  Finally, we emphasize that the main contribution of this paper lies in the
  theoretical analysis of weighted regularization for Fourier extensions and in
  the demonstration of derivative convergence.
  Issues related to fast algorithms and large-scale implementations are largely
  orthogonal to the theoretical developments presented here and will be addressed
  in future work.
\end{itemize}


\section*{Funding}
The research is partially supported by National Natural Science Foundation of China (Nos.
12261131494, 12171455, RSF-NSFC 23-41-00002).
\section*{Author information}
 \subsection*{ \bf Authors and Affiliations}
 School of Mathematics and Statistics, Shandong University of Technology, Zibo, 255049, China \\
 Zhenyu Zhao \& Xusheng Li\\
 Key Laboratory of Deep Petroleum Intelligent Exploration and Development,
Institute of Geology and Geophysics, Chinese Academy of Sciences, Beijing, 100029,China\\
Yanfei Wang\\
Department of Mathematics, Faculty of Physics, Lomonosov Moscow State
University, Vorobyevy Gory, 119991 Moscow, Russia\\
 Anatoly G. Yagola
 \subsection*{\bf Corresponding author}
 Correspondence to Yanfei Wang
 \section*{Ethics declarations}
 \subsection*{\bf Conflict of interest}
 The authors declare that they have no Conflict of interest.
 \subsection*{\bf Non-financial interests}
 The authors have no relevant financial or non-financial interests to disclose.
\end{document}